\documentclass[12pt]{article}
\usepackage{amsmath}
\usepackage{amssymb}
\usepackage{dsfont} 
\usepackage[normalem]{ulem}

\usepackage{graphicx}

\newcommand{\mytopset}{0pt}\newcommand{\myitemsep}{0pt}\newcommand{\myparsep}{0pt}\newcommand{\mypartopsep}{0pt}
\usepackage{enumitem}  

\newcommand{\zparticle}{ZRP-particle} \newcommand{\zparticles}{ZRP-particles}

\newcommand{\mnogitel}{\tau}

\usepackage[dvipsnames]{xcolor} \usepackage{color}

\newcommand{\be}{\begin{equation}} 
\newcommand{\ee}{\end{equation}}
\newcommand{\bel}[1]{\begin{equation}\label{#1}} 
\newcommand{\bea}{\begin{eqnarray}} 
\newcommand{\eea}{\end{eqnarray}}
\newcommand{\balign}{\begin{align}}
\newcommand{\ealign}{\end{align}}
\newcommand{\ba}{\begin{array}}
\newcommand{\ea}{\end{array}}
\newcommand{\bfig}{\begin{figure}}
\newcommand{\efig}{\end{figure}}

\newcommand{\eref}[1]{(\ref{#1})}

\allowdisplaybreaks

\newcommand{\exval}[1]{\mbox{$\langle  {#1} \rangle$}}

\newcommand{\bfP}{{\mathbf P}}

\newcommand{\indic}{1\!\!{\rm I}}   
\newcommand{\bfx}{\mathbf{x}} 
 \newcommand{\bfy}{\mathbf{y}}
\newcommand{\bfz}{\mathbf{z}} \newcommand{\bfn}{\mathbf{n}}
 
\newcommand{\floor}[1]{\lfloor{#1}\rfloor} 

\newcommand{\rmd}{\mathrm{d}}
\newcommand{\rme}{\mathrm{e}}

\newcommand{\ddt}{\frac{\rmd}{\rmd t}}

\newcommand{\ddu}{\frac{\rmd}{\rmd u}}

\newcommand{\R}{{\mathbb R}}

\newcommand{\Z}{{\mathbb Z}}
\newcommand{\N}{{\mathbb N}}
\newcommand{\T}{{\mathbb T}} 

\newcommand{\bzeta}{\boldsymbol{\zeta}}
\newcommand{\bfeta}{\boldsymbol{\eta}}

\newtheorem{theo}{Theorem}[section]
\newtheorem{lmm}[theo]{Lemma}
\newtheorem{df}[theo]{Definition}
\newtheorem{prop}[theo]{Proposition}
\newtheorem{cor}[theo]{Corollary}
\newtheorem{rem}[theo]{Remark}

\newcommand{\proof}{\noindent {\it Proof: }}
\newcommand{\qed}{\hfill$\Box$}

\newcommand{\partnum}{\mathcal{N}}
\newcounter{rnb} 
\newcounter{gnb} 

\begin{document}

\title{Asymmetric exclusion process with long-range interactions}
\author{
Vladimir Belitsky$^{\ast}$ \and 
Ngo P.N. Ngoc$^{\dagger,\ddagger}$ \and 
Gunter M.~Sch\"utz$^{\S}$ 
}

\maketitle

{\small
\noindent $^\ast$ Instituto de Matem\'atica e Est\'atistica,
Universidade de S\~ao Paulo, Rua do Mat\~ao, 1010, CEP 05508-090,
S\~ao Paulo - SP, Brazil, email: belitsky@ime.usp.br\\[1mm]
\noindent $\dagger$ Institute of Research and Development, Duy Tan University, Da Nang 550000, Vietnam\\[1mm]
\noindent $\ddagger$ Faculty of Natural Sciences, Duy Tan University, Da Nang 550000, Vietnam, 
email: ngopnguyenngoc@duytan.edu.vn\\[1mm]
\noindent $\S$ 
Centro de An\'alise Matem\'atica, Geometria e Sistemas Din\^amicos, Departamento de Matem\'atica, Instituto Superior T\'ecnico,
Universidade de Lisboa,
Av. Rovisco Pais 1,
1049-001 Lisbon,
Portugal, email: gunter.schuetz@tecnico.ulisboa.pt}
%

\begin{abstract}
We introduce the headway exclusion process
which is an exclusion process with $N$ 
particles on the one-dimensional discrete torus with $L$ sites
with 
jump rates that depend only on the distance to the next 
particle in the direction of the jump but not on $N$ and $L$. For measures with a long-range 
two-body interaction potential that depends only on the distance 
between neighboring particles we prove a relation between the 
interaction potential and particle jump rates  that is necessary 
and sufficient for the measure to be invariant for the process. 
The normalization of the measure and the stationary
current are computed both for finite $L$ and $N$ and in the 
thermodynamic limit. For a finitely many particles that evolve on 
$\Z$ unidirectionally
it is proved by reverse duality that a certain family 
of non-stationary measures with a microscopic shock and antishock 
evolves into a convex combination of such measures with weights
given by random walk transition probabilities. On macroscopic 
scale this domain random walk is a travelling wave phenomenon 
tantamount to phase separation with a stable shock and a stable 
antishock. Various potential applications of this result and open 
questions are outlined.\\[.3cm]
\textbf{Keywords:} Asymmetric exclusion processes, long-range 
interactions, duality, antishocks, 
phase separation  \\[.3cm]
\textbf{MSC 2020 subject classifications:} Primary 60K35 Secondary: 60K40, 
82C05, 82C22
\end{abstract}

\newpage

\section{Introduction}
\label{Sec:Intro} 

The classical Simple Exclusion Process (SEP) \cite{Spit70,Ligg99} defined on the integer 
lattice $\Z$ or on the finite integer torus $\mathbb{T}_L=\mathbb{Z}/(L\mathbb{Z})$ for some integer 
$L$ is a Markov process where identical particles move on the sites of the lattice and 
satisfy the exclusion rule that allows for no more than one particle per site. Each particle 
attempts to jump at time-points of its Poisson Point Process with the intensity $(r+\ell)$ 
where $r, \ell$ are parameters that determine a jump bias. Independently of everything else, 
a particle attempts  to jump to the neigboring lattice site on its left (resp., right)
with the probability $\ell/(r+\ell)$ ($resp., r/(r+\ell)$). 
The attempt is successful if that target site is not occupied
by another particle, otherwise the jump attempt fails and the particle stays where it was.
This makes the SEP a continuous-time Markov
chain where the total number of particles on the lattice is conserved 
and which is nonreversible if $r\neq \ell$. 

As a model for a physical interacting particle system the exclusion rule
corresponds to an excluded volume interaction (also called hard-core repulsion) 
between particles and 
the bias in the jump rates that arises for $r\neq \ell$ mimicks a driving force that acts
homogeneously on all particles and keeps the particle system
permanently out of thermal equilibrium. 
A very large number of generalizations of this basic process have been
studied. These include driven lattice gas models on $\Z$ or $\T_L$
that allow for particle jumps that depend not only on the occupation of the
target site, but also on the location of other particles. A
prominent example is the Katz-Lebowitz-Spohn model
with next-nearest-neighbour interactions \cite{Katz84}.

Indeed, most investigations have focussed on interactions with finite range,
i.e., where the rate of jump depends only on the particle locations
in a fixed region around the particle that attempts to jump. 
Here we retain exclusion of particles and nearest-neighbor jumps
but lift the 
restriction on the interaction range: We allow the jump rates to depend 
in a general fashion on the distance to the nearest particle in the 
direction of the jump, no matter how far this nearest particle is.
It is convenient to express the distance in terms of the headways, i.e., the number of vacant sites between two particles,
and for this reason we call this process the headway exclusion
process (HEP), see the Markov generators \eqref{AHEPgen} and 
\eqref{genZnew} below. Since we have in mind screened physical interactions only between nearest-neighbor particles
we require the rates of the particle jumps to be independent of both the total
number $N$ of particles and the length $L$ of the torus. 
The latter requirement simply means that 
the microscopic interactions between the particles are pure two-body
interactions that do not ``feel'' the size of the container (the torus in the present case) 
to which they are confined.
It is then easy to see that the headways evolve according
to the zero-range process (ZRP) \cite{Spit70,Andj82,Kipn86}, and various results presented in this work
for the HEP can be proved using known results for the
invariant measure of the ZRP \cite{Spit70,Andj82,Gros03,Evan05,Levi05,Arme09}.

Despite this straightforward mapping between HEP and ZRP it is interesting to study the HEP itself, both from a mathematical perspective and its
potential for applications, some of which are discussed biefly in the last section
of this work. Since the choice of the asymmetry parameter
$r$ and $\ell$ is important for the properties of the process
we say symmetric headway exclusion
process (SHEP) when $r=\ell$, 
asymmetric headway exclusion
process (AHEP) when $r\neq\ell$ and $r\ell \neq 0$, 
and totally asymmetric headway exclusion
process (TAHEP) when  $r\ell = 0$.
Informally, the main results are the following.

\noindent \underline{(i) Invariant measure and stationary current}:
Given the jump rates of the HEP it is easy
to guess from the correspondendence to the ZRP 
the invariant measures of the HEP that we call 
headway measures. It is not obvious, however, whether there 
are other
headway exclusion processes for which the same headway measure is invariant. Theorem~\ref{Theo:invariance} in 
Sec.~\ref{Sec:invmeasure} establishes that the parameters of the
headway measure uniquely define the jump rates of the HEP under the assumption
that these jump rates do not depend on particle number $N$ and
system size $L$ up to the jump asymmetry and a trivial multiplicative
factor that sets the time-scale of the process.\footnote{In fact, the proof also shows that there other HEPs with rates that
depend explicitly on $L$ and $N$, an observation that may be interesting for the study of generalizations of the ZRP but which we do not explore further in the present work.}

From the invariant measure we derive the stationary particle current, see Proposition \ref{Prop:jstatLN}
using a Palm measure approach for fixed $N$ and $L$,
and,  by equivalence of ensembles for the ZRP,  for the thermodynamic limit at 
constant density $\rho=N/L$ Theorem (\ref{Theo:current}). 
It turns out that the phenomenon of condensation that is known from the 
ZRP has a counterpart in the invariant measure of the HEP where it
appears in a somewhat counterintuitive property of the stationary current: For sufficiently
{\it strong long-ranged interaction} potentials (e.g. logarithmic in the distance)
the stationary current becomes linear as a function of the density below some critical particle density $\rho_c$ as if the particles
were {\it non-interacting}.

\noindent \underline{(ii) Random walking density domain:}
To explore this counter intuitive phenomenon in some more detail we switch
gears in Sec.~\ref{Sec:duality}. We focus on the totally asymmetric HEP with $\ell=0$ defined on $\Z$ with an arbitrary
but finite number of particles. We employ a recently developed generalization
of duality \cite{Schu23} to prove on microscopic level a travelling wave 
phenomenon that on macroscopic scale is equivalent to phase separation: 
Theorem \ref{Theo:PS} establishes that a family of measures that is concentrated
on a finite domain on $\Z$ evolves at microscopic time $\tau$ into a convex combination
of such measures with weights given by the transition probabilities at time $\tau$ of a
simple (one-particle) totally asymmetric random walk that moves with a mean velocity $v=r$. Thus one gains full control
over the time-dependent probability measure of the TAHEP
and obtains access to all equal-time correlations.

One can also infer precise information on the hydrodynamic 
limit which has not been established rigorously yet for the ZRP
in the condensation regime. The domain random walk
implies 
for Riemann initial conditions with initial density $\rho(0,x)
= \rho_c \mathbf{1}_{[a,b]}(x)$ that under hyperbolic scaling 
the density at macroscopic time $t$
is the translation $\rho(t,x)
= \rho_c \mathbf{1}_{[a+vt,b+vt]}(x)$ of the initial density profile.
The density discontinuity at the left edge $x_-(t) = a+vt$ 
of the travelling domain constitutes a shock
discontinuity satisfying the Rankine-Huginot condition and is expected
from standard arguments for hyperbolic conservation laws \cite{Lax73,Serr99,Bres02}. On the other hand, the density discontinuity 
at the right edge $x_+(t) = b+vt$  is a stable antishock which is 
surprising since such an discontinuity in the initial state would generically be expected to evolve into a smooth
rarefaction wave \cite{Kipn99}.

The paper is organized as follows. In Sec.~\ref{Sec:AHEP} we first introduce
the basics necessary to define the HEP, viz., lattice and state space along
with the necessary terminology and notation, and then define the HEP both informally and
in terms of its Markov generator. In Sec.~\ref{Sec:invmeasure} we first introduce
the headway measures and state and prove the main theorems concerning
the relationship between the jump rates and the interaction protential for the general HEP
and the properties of the stationary current in the AHEP. In Sec.~\ref{Sec:duality}
we prove for the TAHEP defined on $\Z$ a reverse duality between measures that has the interpretation of a
random walking density domain. Various applications and open questions
are mentioned briefly in Sec.~\ref{Sec:conc}. For self-containedness, the definition
of the zero range process and known properties of its invariant measures
are summarized in the Appendix \ref{App:ZRP}.

\section{The HEP}
\label{Sec:AHEP}

To fix special notation and avoid sending the reader to other publications
for for widely but not universally used conventions we explain
the notation used in the present work in some detail. 
\subsection{Notations, auxiliary definitions and conventions}
\label{Sec:configs}
\subsubsection{Sets of numbers}
\label{Sec:numbers}
The set of strictly positive integers is denoted by $\N$,  and $\N_0 :=\N \cup \{0\}$ 
denotes the set of non-negative integers. The symbol $\N_{m,n}:=\{m,m+1,\dots,n\} \subset \N_0$ 
denotes the integer interval that extends from $m\in\N_0$ to $n\in\N_0$. $\R^+  \equiv (0,\infty)$ 
denotes the set of strictly positive real numbers and $\R^+_0 := \R^+ \cup \{0\} \equiv [0,\infty)$.

\subsubsection{The lattice}
For each $L \in \N$, we denote by $\T_L:=\mathbb{Z}/L\mathbb{Z}$ the discrete 
torus in $\mathbb{Z}$ of size $L$, i.e.,
the one-dimensional lattice 
with $L$ sites and periodic boundary conditions. 
The integers $0,1,\ldots, L-1$ are
called {\it sites}. The sites $0$ and $L-1$ as well as any pair of sites $x$, $y$ such that 
$|y-x|=1$ are declared {\it neighbors}. The site $(x+1)\bmod{L}$ is called the 
{\it right neighbor} of site $x$ and the site $(x-1)\bmod{L}$ is its 
{\it left neighbor}. 

We also consider the infinite integer lattice  $\Z$. In analogy to the torus, integers are called 
{\it sites} and any two sites such that $|y-x|=1$ are declared to be neighbors,
with $x+1$ and $x-1$ being the right and the left neighbor
respectively of 
$x\in \Z$. 

In what follows, $\Lambda$ will be a general notation that means either $\T_L$ or the integer lattice $\Z$.

\subsubsection{Particle configurations}\label{Subsub:particles}
Identical particles are placed on the sites of a lattice $\Lambda$ in such a way 
that each site can be occupied by at most one particle. This constraint (called 
{\it exclusion rule}) means that a configuration of particles on a lattice (to be called 
{\it configuration} and denoted by bold letters like $\bfeta$) can be defined by {\it occupation 
numbers} $\eta_{x}\in\{0,1\}$ where the value $\eta_{x}=1$ is interpreted as presence 
of a particle at site $x\in\Lambda$ and the site is said to be {\it occupied}, while for 
$\eta_{x}=0$ site $x$ is interpreted as {\it vacant}, or {\it empty}, or said to be {\it occupied 
by a vacancy}. 
For the torus, i.e., $\Lambda=\T_L$, 
a configuration is specified by 
the ordered $L$-tuple $\bfeta = (\eta_{0},\dots,\eta_{L-1})$. We call this presentation of 
The {\it particle number} 
is the sum
\begin{equation}\label{Ndef} 
\partnum(\bfeta) := \sum_{x \in\Lambda} \eta_{x} 
\end{equation} 
of occupation numbers. 

Alternatively, a configuration $\bfeta$ can be specified by the set of all sites 
$x\in\Lambda$ such that $\eta_{x}=1$. We call this specification the {\it coordinate presentation}. 
Its formal construction uses a numbering/labelling of particles that is defined as follows. 
Given arbitrary $L$ and $\bfeta=(\eta_{0},\dots,\eta_{L-1})$, we look at site $0$ 
of $\mathbb{T}_L$. If $\bfeta$ has a particle at this site, we attribute to it 
number $1$. Otherwise, we move to the right along $\mathbb{T}_L$ (the rightward direction is 
determined by the position of the right neighbor site) and attribute 
number $1$ to the first particle that we meet. After number $1$ is attributed, we 
continue moving to the right and enumerate sequentially all other particles of 
$\bfeta$. A configuration with $N$ particles is then denoted by the $N$-tuple 
$\bfx=(x_1,x_2,\dots, x_N)$ of particle locations with the $x_i\in\T_L$ in 
increasing order and $x_1\geq 0$ for the position of the ``leftmost'' particle 
and $x_N\leq L-1$ for the position of the ``rightmost'' 
particle.\footnote{It should be borne in mind that the assignments ``leftmost'' and 
``rightmost'' are somewhat artificial for the torus. They are {\it not} meant 
as tags for any specific particle but merely refer to the particle in the 
configuration $\bfx$ with the smallest (largest) coordinate value in the set 
$\N_{0,L-1}$ of coordinate values $x_i\in\T_L$. With this terminology the 
rightmost particle becomes the leftmost after a jump from site $L-1$ to site $0$.}

The coordinate presentation is also used for the configurations with a finite 
number of particles for $\Lambda = \Z$. In this presentation, $x_1\in\Z$ is the position of
the leftmost particle and $x_N\in \Z$ of the rightmost one, where $N=\partnum(\bfeta)$ 
is the particle number.


\medskip On both lattices, 
the relationship between occupation presentation and coordinate presentation 
can be understood as a bijective map and we shall sometimes write $\bfx(\bfeta)$ or 
$\bfeta(\bfx)$ when we change between presentations.
\begin{rem}\label{Rem:Conven} 
Due to the periodicity of the torus, in the presentation by occupation numbers 
of a configuration $\bfeta=(\eta_0, \dots, \eta_{L-1})$ of particles on $\T_L$, the 
indexation is understood modulo $L$, i.e., $\eta_{x+kL} = \eta_{x}$ for any  $k\in \Z$. 
The particle labels
in the coordinate presentation $\bfx=(x_1, \dots, x_N)$ of an $N$-particle configuration are counted modulo $N$, i.e., $x_{i+kN}=x_i$ for any $k\in \Z$.
\end{rem}

\subsubsection{Neighboring particles, headways and distances}\label{Subsub:headways}
Above, the notion of neighboring sites was introduced
to specify
particle configurations on the lattice. 
Here, we define neighborship relations between particles and a related concept called 
{\em headway}. It is used in next sections to construct particle jump rates and the measures 
that are proved to be stationary for the processes studied.  

\paragraph{Integer lattice $\Z$:} 
We call any two particles on $\Z$ {\it neighbors} if there is no 
other particle at any lattice site between them, i.e., in an $N$-particle configuration
$\bfx=(x_1,\dots,x_N)$ with $N\geq 2$, the particles labelled $i$ and $i+1$ 
are neighbors, if $1\leq i \leq N-1$. The leftmost and the rightmost particles (that is,
the particles with labels $1$ and $N$ respectively) have only one 
neighbor, all other particles have two neighbors. If $N=1$, then the only particle
of a configuration has no neighbors.

For $i=1,2, \ldots,N-1$, and any $N$-particle configuration $\bfx=(x_1, \ldots, x_N)$ 
on $\Z$, the particle labelled $i+1$ is called the {\em right neighbor} of the particle labelled $i$, and 
the number of empty lattice sites between $i$-th and $(i+1)$-st particles is called 
{\em headway of the $i$-th particle}; it is denoted by $n_i$ and its formal definition is:
\begin{eqnarray}
n_{i}(\bfx)\, \left(n_i\hbox{ in abbreviated form} \right) \,&: =& x_{i+1} - x_{i} - 1\hbox{ for } 
i\in\N_{1,N-1}, \hbox{ and } \label{hwZdef} \\
n_N(\bfx)\,\left(n_N\hbox{ in abbreviated form }\right)\, &:=&\infty,\label{hwZinfty}
\end{eqnarray} 
where the value $\infty$ for $n_N$ is attributed for  convenience.
We state explicitly that for the configuration
representing the empty lattice  
headways are not defined. However, the headway indicator function introduced below in (\ref{headijdef}) will be defined.

\paragraph{Torus $\T_L$:}
For the torus, these definitions need to be readjusted as explained 
in Remark~\ref{Rem:Conven}. 
In a configuration with $N\geq 2$ particles located at sites $0 \leq x_1 < x_2 < 
\cdots < x_N \leq L-1$ 
any pair of particles with positions $x_{i}$ and 
$x_{i+1}$ as well as the pair with positions $x_{1}$ and $x_{N}$ are 
{\it neighbors}; the particle labelled $(i+1)\bmod{N}$ is the {\em right neighbor}
of the particle labelled $i$ and the particle labelled $(i-1)\bmod{N}$ is its 
{\em left neighbor}. For $N=1$ the (single) particle it its own neighbor and also 
the left and the right neighbor of itself.

For a configuration  $\bfx=(x_1, \ldots, x_N)$ with $N$ particles on $\T_L$, 
{\em the headway of the  $i$-th particle} is the number of the empty sites between it
and its right neighbor (that is the particle labelled $(i+1)\bmod{N}$); formally, it is 
\begin{equation}\label{hwdef}
n_{i}(\bfx)\, \left(n_i\hbox{ in abbreviated form }\right)\, := 
(x_{(i+1)\bmod{N}} - x_{i} - 1) \bmod{L}, \quad i\in \N_{1, N}
\end{equation}
Differently from the lattice $\Z$, now each of $N$ particles possesses its ``natural'' headway and 
they all are finite. Moreover, on $\T_L$ the following relation obviously holds
\begin{equation}\label{distsum}
\sum_{i=1}^N n_i=L-N  \hbox{ for any $N$-particle configuration s.t. }N\leq L. 
\end{equation}
It will be used several times in our proofs.

For clarity, we state explicitly that {\em (a)} \, if there is only one particle (i.e., $N=1$) on $\T_L$ then
it has label 1 and $n_1=L-1$, and {\em (b)} \, if $\T_L$ is totally empty (i.e., $N=0$) then the 
particle headways are not defined, but, as already noted, 
the headway indicator function 
will be defined. 

To illustrate these notions we provide an example: if $\bfeta$ has two 
particles (i.e., $N=2$) on $\T_L\equiv\T_5=\{0,1,2,3,4\}$ positioned at the sites $1$ and 
$3$ then $\bfeta=(0,1,0,1,0)$ and $\bfx(\bfeta)=(x_1, x_2)$ with $x_1=1,\, x_2=3$, that is, the 
particle at site $1$ gets the label ``first'' while the particle at the site $3$ gets the label 
``second''; as for the headways, it holds that $n_1=x_2-x_1-1=1$, i.e., one empty lattice site 
separates the first particle from its right neighbor (this neighbor is the second particle), and 
$n_2=(x_{(2+1)\bmod{N}}-x_2-1)\bmod{5}=2$, i.e., two empty lattice sites separate the second particle 
from its right neighbor (this neighbor is the first particle).

\paragraph{Distances.} For an $N$-particles configuration $\bfx=(x_1, \ldots, x_N)$ on 
$\Lambda$ the {\em distance between the particle labelled $i$ and its right neighbor} is denoted 
by $d_i$ and defined as follows:
\begin{equation}\label{Distance}
d_i:=n_i+1, \hbox{ for }i\in \N_{1,N}\hbox{ if }\Lambda=\T_L, \hbox{ and for }i\in \N_{1,N-1}\hbox{ if }\Lambda=\Z. 
\end{equation}
In our further constructions and proofs, both headways and distances are used although for any pair of 
particles that are neighbors these entities are rigidly connected one to another. The reason is that 
headways frequently appear in studies of exclusion processes whereas distances are used 
for the interaction potential in the construction of Ising-like measures. 

\subsubsection{State spaces for exclusion particles}\label{SpaceEP}
The sets that serve as the state spaces of the process that we study are the domains of all the functions that appear in this work, like, 
for example,  occupation number, headway, and distance
introduced above. The elements of the sets, and thus, the arguments 
of the functions, are configurations that we represent either in the occupation form (and denote by bold Greek letters $\bfeta$, $\bzeta$, etc.) or in the 
coordinate form (and denote by bold Roman letters $\bfx$, $\bfy$, etc.). When no confusion can arise, 
we omit the argument of functions.

\paragraph{Torus $\T_L$:}
For any $L\in\N$ and $\Lambda = \T_L$, the finite set
\begin{equation}\label{gcdef}
\Omega^{gc}_{L}:=\{0,1\}^L
\end{equation} 
is called {\it grand canonical state space}. Its cardinality is given by 
$|\Omega^{gc}_L| = 2^L$, corresponding the number of possibilities of 
placing any number $N\in\{0,\dots,L\}$ of identical exclusion particles on $L$ sites.

For $N\in \N_{0,L}$ the {\it canonical state space}
\begin{equation}\label{candef}
\Omega^{can}_{L,N} := \{\bfeta\in \Omega^{gc}_{L}: \partnum(\bfeta) = N\}
\end{equation}
is the set of all configurations $\bfeta$ with $N$ particles on $\T_L$  as defined in (\ref{Ndef}). We state explicitly that 
$\Omega^{can}_{L,0}$ consists of the single configuration that represents the empty lattice and
\[\Omega^{gc}_{L} = \cup_{N=0}^{L} \Omega^{can}_{L,N}. \]
The cardinality of $\Omega^{can}_{L,N}$ is $\binom{L}{N}$, 
corresponding to the number of possibilities of distributing a fixed number $N$ of identical exclusion 
particles on a fixed number of $L$ sites. We also note that in coordinate presentation  
$\Omega^{can}_{L,N}$ is the set of configurations $\bfx$ with $N$ coordinates $x_i$ such that
\[
0 \leq x_1 < x_2 < \dots < x_N \leq L-1.
\]

\paragraph{Integer lattice $\Z$:} The full state space $\Omega:= \{0,1\}^{\Z}$ of 
all particle configurations on $\Z$ is not countable and not considered in this work. We denote by 
\[
\Omega_N^{x^\ast} := \left\{\bfeta\in\Omega: \eta_{x}=0 \, \forall \,x\in\Z{\hbox{ s. t. }}x <x^\ast,\,
\eta_{x^\ast} = 1, \mbox{ and } \sum_{x=x^\ast}^\infty \eta_{x} = N\right\}
\]
the countable set of all configurations with $N\in\N$ particles, the leftmost of which is 
located on site $x^\ast\in \Z$. In coordinate presentation $\Omega_N^{x^\ast}$ is the set of
configurations $\bfx$ with $N$ coordinates $x_i$ such that
\[
x^\ast = x_1 < x_2 < \dots < x_N 
\]
without further constraint on $x_N$.
The trivial state space corresponding to the empty lattice with $\eta_{x}=0$
for all $x\in\Z$ is denoted by $\Omega_0$. We also introduce the state space
\begin{equation}\label{SpaceForZ} 
\Omega^\ast_N:=\cup_{x^\ast\in \Z} \Omega^{x^\ast}_N,\quad N\in \N 
\end{equation}
of all configurations with $N$ particles on $\Z$.

\subsubsection{State spaces for particles without exclusion}\label{SpaceZRP}
%
We define 
\begin{equation}\label{TiOmegaGC}
\tilde{\Omega}^{gc}_{N} := \N_0^N, \quad N\in \N, 
\end{equation}
and the definition implies that an element from $\tilde{\Omega}^{gc}_{N}$ is an $N$-tuple
$\bfn = (n_{0}, \dots , n_{N-1})$ where each $n_i\in \N_0$. We interpret $n_i$
as the number of identical particles that are piled up at the site $i$ of $\T_N$, 
and thus, an element from $\tilde{\Omega}^{gc}_{N}$ can be understood as representing
a configuration of identical particles put on sites $\T_N$ {\it without} exclusion, 
and, accordingly,  we call $\tilde{\Omega}^{gc}_{N}$ {\em the grand canonical state 
space of particles without exclusion}. Then for $K\in\N_0$,  we denote by $\tilde{\Omega}^{can}_{N,K}$ the finite subset of $\tilde{\Omega}^{gc}_{N}$ 
with a total number of $K$ particles without exclusion on $\T_N$. This canonical state 
space is finite and has cardinality $|\tilde{\Omega}^{can}_{N,K}| = {N+K-1 \choose K}$. 

The elements $\tilde{\Omega}^{can}_{N,K}$ can be associated to sets of the headways 
of the configuration with $N$ exclusion particles distributed on a torus
of $L=K+N$ sites. This association is one of the ingredients of the 
proof of Proposition~\ref{Prop:Zcan} of Section~\ref{Sec:partfunc}. Combining this proposition with properties of the zero range process (ZRP) \cite{Spit70,Kipn99} briefly reviewed in App. 
\ref{App:ZRP} yields the 
current-density relation for the TAHEP that is stated in Theorem~\ref{Theo:current}. In the arguments
and constructions that lead to this theorem, it will become clear that ``particles without
exclusion'' mentioned above are actually the ZRP particles.
%
%

\subsection{Definition of the HEP}

\subsubsection{Informal description}
\label{Sec:process}

Since the number of particles does not change over time, 
we fix from now on an arbitrary integer $N$ and
describe informally the stochastic dynamics of $N$ particles. The rules are determined by two parameters 
$r,\ell\in \R_0^+$ and a sequence
\begin{equation}\label{wrates}
w_0=0, w_n \in \R_0, n\in \N  
\end{equation}
of positive parameters that for reasons that will become
clear below are called jump rates.
Let now $t\in\R_0$ stand for an arbitrary positive time and let $\bfx=(x_1, \ldots, x_N)$ denote 
the particle positions on $\T_L$ at time $t$. 
Each particle is equipped with two alarm clocks and all
$2N$ clocks 
are independent of each other and their
alarm times are independent
exponentially distributed random variables. 
For the $i$-th particle, the parameter of clock 1 
depends on the number of empty sites to the right 
of the particle, viz., $r w_{n_i}$, and the parameter of clock 2  
depends on the number of empty sites to the left
of the particle, viz., $\ell w_{n_{i-1}}$.
If clock 1 (clock 2) of particle $i$ rings first, this particle immediately
jumps to the right (left) neighboring site.
It is important to notice that when a headway is $0$ then the 
parameter of the corresponding 
exponential distribution is by definition $w_0=0$.
This means that the corresponding alarm 
will go off at time $\infty$ (i.e., it never goes off). Hence
no jump to an occupied site is ever attempted
which ensures that the exclusion rule is respected
at all times $t\in\R_0$.
The jump rearranges the particle position of particle $i$ and makes the
particles $i-1$, $i$, and $i+1$ readjust instantaneously their alarm clock parameters.

The definition of the TAHEP on $\Z$ is analogous. The differences are: $r=1$, $\ell=0$, and 
the set (\ref{wrates}) is enlarged by adding one more number denoted by $n_{\infty}$ which
is needed because the headway of the rightmost particle is $\infty$. 

\subsubsection{Formal definition}
\label{Sec:generator}
Following the standard procedure of construction of interacting particle systems (see 
Chapters 2 - 4 in \cite{Ligg10} and Chapter~I in \cite{Ligg85}) we define both the
HEP on $\T_L$ and the TAHEP on $\Z$ in terms of the generators of the 
corresponding Markov semigroups. 

We introduce the the auxiliary notation $\bfeta^{x,y}$ that stands for the configuration
obtained from configurations $\bfeta$ of particles on $\Lambda$ via permuting the 
occupation numbers $\eta_{x}$ and $\eta_{y}$, i.e., the configuration 
which occupation number at each $z\in \Lambda$  is as follows:
\[
\eta^{x,y}_z = \eta_{z} + (\eta_{y}-\eta_{x}) 
(\delta_{z,x} - \delta_{z,y})
= \left\{\begin{array}{ll} \eta_z &  z\notin \{x,y\}, \\
\eta_{y} &  z=x, \\
\eta_{x} &   z=y.\end{array}\right. 
\] 
Another auxiliary tool are the {\it headway indicator} functions 
\begin{equation}\label{headijdef}
h_{x,n}(\bfeta) = \left\{
\begin{array}{ll}
1 & \hbox{ iff } \bfeta_x=1\hbox{ and the headway of the particle at site }x\hbox{ is }n \\
0 & \hbox{ else.}
\end{array}
\right.
\end{equation}
When working on $\T_L$ all site indices that appear in these definitions are understood modulo $L$.
When we use $h_{x,n}(\cdot)$ on $\Omega^{gc}_L$, the range of $x$ is $\N_{0,L-1}$ and the range of $n$ 
is $\N_{0,L-1}$ (the longest possible headway is achieved when there is one particle on $\T_L$ and it is 
$L-1$). When we use $h_{x,n}(\cdot)$ on $\Omega^{\ast}_N$, the range of $x$ is $\Z$ and the range of $n$ 
is $\N_{0}\cup\{\infty\}$. It will be always clear from context what the domain of $h_{x,n}(\cdot)$
is, and, consequently, what the ranges of its index variables are. Notice that for the configuration representing 
the  empty lattice the value of $h_{x,n}$ is always $0$ because $\eta_x=0$
for all sites $x$. 

With these definitions we are in a position to define the processes that we study.
 
\begin{df}{(Headway exclusion process on the integer torus).}\label{Def:onTorus}
For any $L\in \N$, any $r,\ell\in\R_0^+$, and any set $\{w_n, n\in \N_0\}$ satisfying (\ref{wrates}),
the continuous-time Markov chain with the state space
$\Omega^{gc}_L$ whose Markov generator $\mathcal{L}$ acts on bounded 
functions $f:\Omega^{gc}_L\to\R$ as
\begin{equation}\label{AHEPgen}
(\mathcal{L} f)(\bfeta) = \sum_{x\in\T_L} \left\{r c^+_{x}(\bfeta)
\left[f(\bfeta^{x,x+1}) - f(\bfeta)\right] +
\ell c^-_{x}(\bfeta)
\left[f(\bfeta^{x,x-1}) - f(\bfeta)\right] \right\} 
\end{equation}
where 
\begin{equation}\label{transrates}
c^{+}_{x}(\bfeta) := \sum_{n=0}^{L-1} w_{n} h_{x,n}(\bfeta), \quad 
c^{-}_{x}(\bfeta) := \sum_{n=0}^{L-1} w_{n} h_{(x-1-n)\bmod L,\,n}(\bfeta)  \quad \forall \bfeta\in \Omega^{gc}_L
\end{equation} 
is called the {\em headway exclusion process (HEP) on $\T_L$}. Specifically, when $r=\ell$ it is called 
symmetric headway exclusion process (SHEP), when
$r \neq \ell$ it is called 
asymmetric headway exclusion process (AHEP),
and when $r=0$ or $\ell=0$ 
it is called 
totally asymmetric headway exclusion process (TAHEP).
\end{df}

\begin{rem}
Since by definition $\bfeta^{x+1,x} = \bfeta^{x,x+1}$
the {\em bond transition rate} 
\begin{equation}
c(\bfeta) := r c^+_{x}(\bfeta) + 
\ell c^-_{(x+1)\bmod L}(\bfeta)
\end{equation}
allows for writing the generator in the form
\begin{equation}\label{AHEPgenalt}
(\mathcal{L} f)(\bfeta) = \sum_{x\in\T_L}  c_{x}(\bfeta)
\left[f(\bfeta^{x,x+1}) - f(\bfeta)\right]  \, \forall \bfeta\in \Omega^{gc}_L.
\end{equation}
The transition rate from a configuration 
$\bfeta$ to a configuration 
$\bzeta$ is then given by
\begin{equation}\label{tarrate}
c(\bfeta\to \bzeta):= 
 \sum_{x\in\T_L} c(\bfeta) \delta_{\bzeta,\bfeta^{x,x+1}}
\left[f(\bzeta) - f(\bfeta)\right] 
\end{equation}
 for any $\bfeta,\bzeta\in \Omega^{gc}_L, \,\bfeta\not=\bzeta$ and can be expressed as
\begin{equation}\label{tarrate2}
c(\bfeta\to \bzeta)=\left\{\begin{array}{l} r w_n \hbox{ \ if and only if }\eta_x=1,\, 
\eta_{(x+1)\bmod L}=0,\hbox{ there are }\\
\hbox{ \ \ \ $n$ empty sites between the particle at $x$ and its right }\\
\hbox{ \ \ \  neighboring particle in $\bfeta$, and } \bzeta=\bfeta^{x, (x+1)\bmod L},\\
\ell w_{n}  \hbox{ \ if and only if }\eta_x=1,\, \eta_{(x-1)\bmod L}=0,\, \hbox{ there are }\\
\hbox{ \ \ \ $n$ empty sites between the particle at $x$ and its left}\\
\hbox{ \ \ \ neighboring particle in $\bfeta$, and } \bzeta=\bfeta^{x, (x-1)\bmod L},\\
0\hbox{ \ \ \ \ \ \ \ \ \ }otherwise.
\end{array}\right.
\end{equation}
\end{rem}

For the TAHEP considered in this work we choose $r=1$
and $\ell=0$ and write $c_x(\bfeta)$ instead of $c^{+}_{x}(\bfeta)$.

\begin{df} {(TAHEP with finite particle number on $\Z$).}\label{Def:onZ}
For any $N\in \N$ and any set $\{w_i, i \in \N_0\cup\{\infty\}\}$ satisfying (\ref{wrates}),
the continuous-time Markov process with the state space
$\Omega^{\ast}_N$ whose Markov generator $\mathcal{L}$ acts on bounded 
functions $f:\Omega^{\ast}_N\to\R$ as
\begin{equation}\label{genZnew}
(\mathcal{L} f)(\bfeta) = \sum_{x\in\Z} c_{x}(\bfeta)
\left[f(\bfeta^{x,x+1}) - f(\bfeta)\right] 
\end{equation}
where 
\begin{equation}\label{ratesZnew}
c_{x}(\bfeta) := \sum_{n=0}^{\infty} w_{n} h_{x,n}(\bfeta)\,+\, w_{\infty} h_{x,\infty}(\bfeta)
\end{equation} 
is called {\em TAHEP with $N$ particles on $\Z$}, or, simply, {\em TAHEP on $\Z$}. 
\end{df}


%

\noindent
{\bf Comment (a)}. The sums in (\ref{transrates}) and (\ref{ratesZnew}) may be taken starting from $1$
rather than from $0$ because $w_0=0$.

\medskip\noindent{\bf Comment (b)}. The headway $h_{(x-1-n)\bmod L,\,n}(\bfeta)$ used in (\ref{transrates})
 is equal to $1$ if and only if there is a particle at site $x$ and 
the number of empty site between it and its nearest-neighboring particle on its {\it left} is $n$. This event 
depends upon the number of the empty sites that a particle at $x$ sees on its left and therefore cannot be 
expressed by $h_{x, m}$ since this function ``looks'' to the right of the particle that occupies $x$.
This discrepancy explains why we put  $(x-1-n)\bmod L$ in the subscript of $h$ that we used.

%
 
\medskip\noindent{\bf Comment (c)}. 
It is straightforward to see 
that for fixed $L\in\N$, the HEP described above is, 
in fact, a collection of independent processes, each one being characterized by the conserved 
particle number $N\in\N_0$ and with irreducible state space $\Omega^{can}_{L,N}$. 
(For $N=0$ and any $L$, this is  the trivial process where nothing happens.) This explains  
why we have defined these spaces (see (\ref{candef})) and why we shall construct and analyse 
measures concentrated on these sets (see Def.~\ref{Def:CanHead}).

\section{Invariant measures of the HEP on $\T_L$}
\label{Sec:invmeasure}
The contents and the structure of the present section are as follows. In Section~\ref{Sec:measure} 
we define a family of probability measures that we call {\it headway measures}; the construction 
employs ideas from thermodynamics. In Section~\ref{Sec:invar} we establish that these measures 
are invariant for the HEP  under certain conditions on its jump rates. Then, in Section~\ref{Sec:partfunc} we turn our 
attention to the normalizing constant of the
headway measures which we express in terms of parameters of measures 
that are invariant for the zero range process. This expression and certain properties of the zero range process 
(that are summarized in Appendix A) allow us then to find a formula for the density-current relation 
in the AHEP when it evolves under a stationary measure, 
see Section~\ref{Sec:current}.

\subsection{Headway measures for particle configurations on $\T_L$}
\label{Sec:measure}
It is convenient to use established terminology borrowed from thermodynamics even though 
the fundamental thermodynamic notion of thermal equilibrium plays no role in the 
present context. Therefore we will ignore the thermodynamic parametrization 
$\beta=1/(k_B T)$ for the temperature $T$ and Boltzmann constant $k_B$ and we will 
set $\beta=1$ in all expressions in which this quantity would appear in thermodynamics.

Let 
\begin{equation}\label{Jpotential}
J(0):=\infty \hbox{ and } J(n)\in\R\cup\{\infty\}, \,\, n\in\N, 
\end{equation}
be an arbitrarily fixed sequence and, for $L\in \N$ define 
\begin{equation}\label{Edef}
E_{L}(\bfeta) := \sum_{x=0}^{L-1} \sum_{n=0}^{L-1} J(n+1) h_{x,n}(\bfeta), \quad \forall\bfeta\in\Omega^{gc}_{L},
\end{equation}
where $h_{x,n}(\cdot)$ is the headway indicator function from (\ref{headijdef}) with domain $\Omega^{gc}_L$.  
In thermodynamic interpretation, $E_L(\bfeta)$ is the internal energy of a {\em microstate} 
$\bfeta$ with {\it static interaction potential $J(r)$ between nearest-neighbor particles 
that are at distance} $r$, $r=0,1,\ldots$.  
Accordingly, we shall call the set (\ref{Jpotential}) {\it interaction potential}. 
Recall from (\ref{Distance}) that if the distance between two nearest-neighbor particles 
is $n+1$ then the headway of one of them is $n$. 
 
For given set $\{J(\cdot)$\} and corresponding {\it Boltzmann factors} 
\begin{equation}\label{Bfactors}
y_0:=0\hbox{ and }y_{r} := \exp{(-J(r))} \in \R^+_0, \,r\in \N,
\end{equation}
we introduce the {\it Boltzmann weight} that will be the basic brick of our 
construction of headway measures:
\begin{equation}\label{BWdef}
\pi_{L}(\bfeta) := \rme^{- E_L(\bfeta)} \equiv \sum_{x=0}^{L-1} \sum_{n=0}^{L-1} J(n+1) h_{x,n}(\bfeta) 
= \prod_{x=0}^{L-1} \prod_{n=0}^{L-1} y_{n+1}^{h_{x,n}(\bfeta)}, \quad \forall\bfeta\in\Omega^{gc}_{L} .
\end{equation}

\begin{df}[Canonical headway measure]\label{Def:CanHead}
For any $L\in \N$ and any $N\in \N_{0, L}$, the measure $\mu^{can}_{L,N}$ defined by 
\begin{equation}
\mu^{can}_{L,N}(\bfeta)  := \frac{1}{Z^{can}_{L,N}} \pi_{L}(\bfeta)
, \qquad \forall \, \bfeta \in \Omega^{can}_{L,N}
\label{headwaycandef}
\end{equation}
is called the {\em canonical headway measure} and the normalization 
\begin{equation}
Z^{can}_{L,N} :=  
\sum_{\bfeta\in\Omega^{can}_{L,N}} \pi_{L}(\bfeta)
= \sum_{\bfeta\in\Omega^{can}_{L,N}} \prod_{x=0}^{L-1} \prod_{n=0}^{L-N} 
y_{n+1}^{h_{x,n}(\bfeta)}
\label{Zcandef}
\end{equation}
is called the {\em canonical partition function}. The {\em canonical ensemble} is
the probability space $(\Omega^{can}_{L,N}, \mathcal{F}^{can}_{L,N},\mu^{can}_{L,N})$, 
where $\mathcal{F}^{can}_{L,N}$ is the $\sigma$-field that contains all the singletons 
of $\Omega^{can}_{L,N}$. 
\end{df}

\begin{df}[Grand canonical headway measure]\label{Def:GranCanHead}
For any $L\in \N$ and any $z\in \R^+_0$, the {\em grand canonical headway measure} $\mu^{gc}_{L,z}$
is defined by 
\begin{equation}
\mu^{gc}_{L,z}(\bfeta) := \frac{1}{Z^{gc}_{L,z}} z^{\partnum(\bfeta)} \pi_{L}(\bfeta)
, \qquad \forall \, \bfeta \in \Omega^{gc}_{L}
\label{headwaygcdef}
\end{equation}
with particle number function (\ref{Ndef})
where $z$ is called {\em fugacity} and the normalization 
\begin{equation}
Z^{gc}_{L,z} := \sum_{\bfeta\in\Omega^{gc}_{L}} z^{\partnum(\bfeta)} \pi_{L}(\bfeta)
\label{Zgcdef}
\end{equation} 
is called  {\em grand canonical partition function}. The {\em grand canonical ensemble} is the probability space $(\Omega^{gc}_{L}, \mathcal{F}^{gc}_{L}, \mu^{gc}_{L,z})$, where 
$\mathcal{F}^{gc}_{L}$ is the $\sigma$-field that contains all the singletons of $\Omega^{gc}_{L}$.
\end{df}

 \begin{rem}\label{Rem:reflection}
Since the Boltzmann weight depends only on the distances
(\ref{Distance}) between pairs of nearest-neighbor 
particles,
it is invariant under space reflection $\eta_{x} \mapsto \eta_{L-1-x}$. Therefore, also the 
canonical and the grandcanonical headway measures are invariant under space reflection. This fact 
will be used to simplify the proof of Theorem~\ref{Theo:invariance}. 
\end{rem}

\begin{rem}\label{Rem:alternative}
For each $L\in \N$ and each $N\in \N_{0,L}$, one can define 
{\em the embedded canonical measure} by (below, $\delta$ is the Kronecker-$\delta$)
\[  \hat{\mu}^{can}_{L,N}(\bfeta)  := \frac{1}{Z^{can}_{L,N}} \pi_{L}(\bfeta) 
    \delta_{\partnum(\bfeta), N}, \quad \forall \, \bfeta \in \Omega^{gc}_{L}  \]
where a configuration $\bfeta\in \Omega^{gc}_{L}$ with $\partnum(\bfeta) \neq N$ particles 
has zero probability in $\hat{\mu}^{can}_{L,N}(\bfeta)$. Our remark is that the grand canonical 
measure may be re-written as the convex combination
\[   \mu^{gc}_{L,z} = \sum_{N=0}^{L} \gamma^{gc}_{L,N}(z) \, \hat{ \mu}^{can}_{L,N}  \]
 of embedded canonical measures with the weight
\[   \gamma^{gc}_{L,N}(z) = \frac{Z^{can}_{L,N}}{Z^{gc}_{L,z}} z^{N} \]
and grand canonical partition function expressed as
\[   Z^{gc}_{L,z} = \sum_{N=0}^{L} \sum_{\bfeta\in\Omega^{can}_{L,N}}
     z^{\partnum(\bfeta)} \pi_{L}(\bfeta) = \sum_{N=0}^{L} z^{N} Z^{can}_{L,N}. \] 
If $\bfeta\in\Omega^{can}_{L, N}$ then the embedded canonical measure $\hat{\mu}^{can}_{L,N}$ 
and the canonical measure $\mu^{can}_{L,N}$ attribute to $\bfeta$ the same weight, while if 
$\bfeta\not \in\Omega^{can}_{L, N}$ then the former gives value $0$ while the latter is not 
defined. This advantage of $\hat{\mu}^{can}_{L,N}$ over $\mu^{can}_{L,N}$ is useful for certain 
situations, one of which is the argument that we shall use at the beginning of the proof of 
Theorem~\ref{Theo:invariance}.
\end{rem}

\subsection{Invariance and uniqueness of the headway measures and uniqueness of the jump rates}
\label{Sec:invar}
The purpose of this section is to prove the following two assertions:\\
\begin{description}[topsep=\mytopset,itemsep=\myitemsep, parsep=\myparsep,partopsep=\mypartopsep]
\item{(P1)} For a given interaction potential (\ref{Jpotential}), a given system size $L$ and a given particle number $N$,  
the canonical headway measure is the unique invariant measure for the 
HEP on $\T_L$ involving $N$ particles with a certain choice of jump rates (\ref{wrates}) 
that will be stated explicitly below.
\item{(P2)} This choice of rates is unique in the sense that it is only
possible choice of rates that does not depend on $L$ and $N$.
The latter statement does not rule out different rates $w_r(L,N)$ for some
fixed values of $L$ and $N$, but as discussed in the introduction and re-enforced in 
Remark~\ref{Rem:OtherMeasures} we are not interested in such rates.
\end{description}

Assertion (P1) follows from the known invariant measure of the ZRP \cite{Andj82}. Hence the proof is omitted. Assertion (P2) 
uses
that 
a measure $\mu$ is time-invariant for a process 
with transition rates $c(x\rightarrow y)$ if and 
only if it satisfies the general balance relation
\begin{equation}\label{wepaL}
\sum_{z\in \Omega} \mu(z)\,c(y\rightarrow z) =
 \sum_{x\in \Omega} \mu(x)\,c(x\rightarrow y), \qquad \forall y\in \Omega,
\end{equation} 
see, for example, Proposition~2.13 in Liggett \cite{Ligg85}.

\begin{theo}
\label{Theo:invariance} 
{\bf (a)} Fix arbitrarily an interaction potential (\ref{Jpotential}) satisfying 
\begin{equation}\label{NotInf}
J(k)\not=\infty,\,\, k\in \N.
\end{equation} 
Let then $\{w_k,\, k\in \N_0\}$ be a set of real numbers that is given by this potential via
the relation 
\begin{equation} \label{rates}
w_{n} = \left\{\ba{ll} 0 & n=0 ;\\   w \, \rme^{J(n+1) - J(n)} & n\in\N ,  \ea\right.
\end{equation}
where $w$ is an arbitrary strictly positive real number. Let next $L\in\N$ and $N\in \N_{0, L}$ be arbitrary and consider 
the HEP with $N$ particles on $\T_L$ and with particle  jump rates (\ref{rates}).
Then, for any value of the fugacity $z$, the grand canonical headway measure $\mu^{gc}_{L,z}$ 
corresponding to the fixed potential $J(k)$ is invariant for the HEP. \\

\noindent {\bf (b)} \ If $\mu^{gc}_{L, z}$ corresponding to some interaction potential (\ref{Jpotential}) 
satisfying (\ref{NotInf})
and some arbitrary value of the fugacity $z$ is invariant for the HEP for every lattice size $L$ 
and every particle number $N\in \N_{0,L}$, then there exists $w\in\R^+$ such that the particle jump 
rates $\{w_k,\, k\in \N_0\}$ of the HEP is related to the potential via (\ref{rates}).
\end{theo}


\noindent{\it Proof.} \ 
First we note that it is sufficient to prove the theorem for the TAHEP. The justification is as follows:
If a headway measure is invariant for TAHEP, then, due to the space reflection symmetry of headway 
measures (see Remark~\ref{Rem:reflection}), the same measure will be invariant for 
the ``reversed TAHEP'', i.e., TAHEP in which particles jump to the left rather than to the right 
(corresponding to the choice $\ell=1, r=0$ in (\ref{tarrate})). 
Since according to Definition~\ref{Def:onTorus} the generator of the HEP is a linear combination of generators 
of the TAHEP and of the reversed TAHEP with the 
respective multipliers $r$ and $\ell$, the headway measure is invariant for the HEP as well. 

Second, we note that invariance of the canonical 
measure $\mu^{can}_{L,N}$ implies invariance of the
grand canonical measure. This is because invariance of the canonical measure manifestly implies invariance of
the embedded canonical measure $\hat{\mu}^{can}_{L,N}$
defined in Remark~\ref{Rem:alternative} which gives
weight zero to the subspace of $\Omega^{gc}_L$ 
on which the process with $M\neq N$ particles evolves.
Hence, since all embedded canonical measures are 
defined on the same space $\Omega^{gc}_{L}$ any convex combination of these measures 
is invariant for the process. 
In particular, it follows that
that the grand canonical measure $\mu^{gc}_{L,z}$ is 
invariant for the process for any value of the chemical potential $z$.
Therefore the theorem can be proved by focussing on the
TAHEP and the canonical 
measure $\mu^{can}_{L,N}$.

%
\medskip\noindent{\it Proof of (b).}  
We analyse the TAHEP with $N=2$ particles on $\T_L$ for any arbitrarily fixed value of  
$L>2$. (The cases $L=1,2$ may be treated by a direct calculation rather than by our generic analysis
because there may only be one or two particles on the lattice.) 
It will transpire that one can choose rates that depend on $L$ for 
$\mu^{can}_{L,2}$ to be invariant, but that the only choice of rates that does not depend 
on $L$ is \eref{rates}. For $N=3$ one can conceivably choose another set of rates. However,
since independence of the rates of $L$ {\it and} $N$ is to be proved, the rates determined by choosing $N=2$ 
determine the rates to be taken for all $N$.

Let $\pi_{L,2}$ denote the non-normalized measure 
 $\mu^{can}_{L, 2}$, and for  a configuration $\bfx=(x_1, x_2)\in \Omega^{can}_{L,2}$ 
of two exclusion particles on $\T_L$, let $\hat{\pi}_{L,2}(r_1,r_2)$ denote the value of $\pi_{L,2}(\bfx)$,
where  $r_1: = d_1=(x_2-x_1) \bmod{L}$ and $r_2: d_2 =(x_1-x_2) \bmod{L}$ are the distances as 
defined in (\ref{Distance}). According to the definition \eref{headwaycandef}, it holds that
\begin{equation}\label{IsingL2}
\hat{\pi}_{L,2}(r_1,r_2) = y_{r_1} y_{r_2} .
\end{equation}
In terms of the distance the rate of jump of particle $i\in\{1,2\}$ is given by 
$w_{r_i-1}$ with $w_{0} = 0$ due to exclusion.

The condition \eref{wepaL} for the invariant measure when applied to the 
configuration $\bfx$, yields the equality
\begin{eqnarray*}
0 & = & w_{r_1} \pi_{L,2}(r_1+1,r_2-1) - w_{r_1-1} \hat{\pi}_{L,2}(r_1,r_2)
\nonumber \\ & &
+ w_{r_2} \pi_{L,2}(r_1-1,r_2+1) - w_{r_2-1} \hat{\pi}_{L,2}(r_1,r_2).
\end{eqnarray*}
Since there are exactly two particles on the torus with $L$ sites, one has $r_1+r_2=L$. Accordingly, 
we can simplify the above expression via the substitution  $r:=r_1$ and $r_2 = L - r$. 
Hence, with \eref{IsingL2}, the invariance of $\pi_{L,2}$ means that the rates must satisfy
\be
0 = w_{r} y_{r+1} y_{L-r-1} - w_{r-1} y_{r} y_{L-r}
+ w_{L-r} y_{r-1} y_{L-r+1} - w_{L-r-1} y_{r} y_{L-r}.
\label{inv2}
\ee
Now set without loss of generality
\[
w_{r} := a_{r} \frac{y_{r}}{y_{r+1}}
\]
with coefficients $a_r$, $r\in\N_{1,L-1}$ and recall that $y_{0} = 0$ 
which enforces $w_{0}=0$ without condition on the coefficient $a_0$. 
After rearranging the negative terms according to the pairwise balance 
relation the stationarity condition for the rates becomes
\[
0  =  (a_{r} - a_{L-r-1} ) y_{r} y_{L-r-1} + (a_{L-r}  - a_{r-1}) y_{r-1} y_{L-r} 
\]
which in terms of the combination $b_r:= (a_{r} - a_{L-1-r}) y_{r} y_{L-r-1} $ reads
\[
0 = b_{r} - b_{r-1}, \quad \forall \, r \in \N_{1,L-1}
\]
Hence
\[
b_{r} = (a_{r} - a_{L-1-r}) y_{1+r-1} y_{1+L-r-2} = c, \quad \forall \, r \in \N_{0,L-1}
\]
with a constant $c$. Taking $r=0$ or $r=L-1$ one finds $c=0$ since $y_{0}=0$ due to exclusion.
Therefore, for any fixed $L>2$ the rates must satisfy the reflection condition
\be 
w_{r} \frac{y_{r+1}}{y_{r}}  
= w_{L-1-r}  \frac{y_{L-r}}{y_{L-r-1}}
, \quad \forall \, 1 \leq r \leq L-2.
\label{refl}
\ee
Specifically, for $L=r+2$ this implies
\[
w_{L-2} \frac{y_{L-1}}{y_{L-2}}  =  w_{1} \frac{y_{2}}{y_{1}}, \quad \forall \, L \geq 3 .
\]
with $w_{1}$ being a free parameter that may be chosen to depend
on $L$. Demanding independence of $L$ and choosing $w_{1}=w 
y_{1}/y_{2}$ with a free parameter $w\in\R^{+}$ fixes all rates $w_r$ to be of the form
\[
w_{0} = y_{0} , \quad w_{r} = w \frac{y_{r}}{y_{r+1}} , \quad r \geq 1, 
\]
in terms of the quantities $y_r=\rme^{-J(r)}$ which is equivalent to
\eref{rates}. \qed

\begin{rem}\label{Rem:OtherMeasures} 
The reflection condition \eref{refl} is a set of $\lfloor L/2\rfloor -1$ equations
does {\it ** not **} require the jump rates for the process with two particles to be 
of the form \eref{rates} to ensures invariance of the
headway measure. This fact demonstrates that other solutions to the invariance 
condition \eref{inv2} may exist. However, requiring independence of $L$ limits the set of solutions 
to rates of the form \eref{rates}.
\end{rem}

\subsection{Partition function of the headway measures on $\T_L$}\label{Sec:partfunc}

The normalization factors are not easy 
to handle in computations when given as in \eref{Zcandef} and \eref{Zgcdef}. It turns out, however, that they are related to the 
normalization factor appearing in the invariant measure of the zero range process. The 
necessary facts about this measure are  described briefly in Appendix A. This relation 
is the essence of Proposition~\ref{Prop:Zcan} stated and proved below. 
It will be then used in Section~\ref{Sec:current} to derive a formula that for
the current-density relation in the AHEP.

Before stating Proposition~\ref{Prop:Zcan} we note in lemma below 
some properties of the headway indicator involving products with occupation 
numbers. Thsee properties are used in the proof of Proposition~\ref{Prop:Zcan} and the arguments
of Section~\ref{Sec:current}. We stress the lemma is stated for $\T_L$, not for $\Z$.

\begin{lmm}
\label{Lem:head} 
Let $L$ be an arbitrary integer and let $\bfeta=(\eta_0, \ldots, \eta_{L-1})$ be an arbitrary 
configuration of exclusion particles on $\T_L$. 
Then the headway indicator  on $\T_L$ defined by (\ref{headijdef}) satisfies 
the following identities 
\begin{eqnarray}
[h_{x,n}(\bfeta)]^m & = & h_{x,n}(\bfeta)\hbox{ for } n\in \N_{0, L-1}, \, x\in \T_L, \, m\in \N, 
\label{head1} \\
\eta_x h_{x,n}(\bfeta) & = & h_{x,n}(\bfeta)\hbox{ for } n\in \N_{0, L-1}, \, x\in \T_L,
\label{head4} \\
\eta_{x+1} h_{x,n}(\bfeta) & = & 0\hbox{ for } n\in \N_{1, L-1},\, x\in \T_L,
\label{head5} \\
\hbox{if }\eta_0=1&\!\!\!\hbox{for}&\! \!\!n\in \N_{0, L-1},\, x\in \T_L\hbox{ s.t. }x\in \N_{L-n, L-1},
\hbox{ then }h_{x,n}(\bfeta)=0.
\label{head3}
\end{eqnarray}
\end{lmm}

\proof The identities follow directly form the definitions of the headway indicator and the occupation number.
Still, however, we present the proofs.

The value of $h_{x,n}(\bfeta)$ is always  either $1$ or $0$. This implies directly (\ref{head1}). Moreover,
for the value to be $1$ it is necessary that $\bfeta$ has a particle at $x$, i.e., $\eta_x=1$. This implies
(\ref{head4}). For the proof of other identities, it is important to notice that for $h_{x,n}(\bfeta)$ to be $1$
it is necessary that $n$ be the headway of the particle at $x$ and that $n>0$. 
If $\eta_{x+1}=1$, then the headway of the particle at $x$ is $0$ and, consequently, $h_{x, n}=0$ for any 
$n\geq 1$, that implies that (\ref{head5}) holds true; if to the contrary $\eta_{x+1}=0$, then 
$\eta_{x+1} h_{x,n}=0$ independently of the value
of $h_{x,n}$, that is, (\ref{head5}) again holds. Finally, to prove (\ref{head3}), one has to notice that 
there are at most $n-1$ sites between the site $x$ and the site $0$ on $\T_L$ (this is because of the requirement 
$L-n \leq x \leq L-1$). Thus, if $0$ contains a particle, then the the headway of the particle at $x$ must 
be smaller than $n$, and, consequently,  $h_{x, n}(\bfeta)=0$. \qed

We are now in a position to present the canonical partition function in terms of the canonical partition function
of the zero range process.

\begin{prop}
\label{Prop:Zcan} 
Fix arbitrarily interaction potential $\{J_i\}$ that satisfies (\ref{Jpotential}). Next, for arbitrary  
$L\in \N$ and  $N\in \N_{0,L}$
consider the canonical headway measure $\mu^{can}_{L,N}$ 
of Definition~\ref{Def:CanHead}) that corresponds to the fixed potential.
Then the canonical partition function  \eref{Zcandef}
satisfies  
\begin{eqnarray}
Z^{can}_{L,0} &=& 1,\nonumber \\
Z^{can}_{L,N} &=& 
\frac{L}{N}\sum_{n_1=0}^{\infty} {y_{n_{1}+1}} \sum_{n_2=0}^{\infty} {y_{n_{2}+1}}
\dots \sum_{n_N=0}^{\infty} {y_{n_{N}+1}} \delta_{\sum_{i=1}^N n_{i},L-N}, \,  N \in \N_{1, L}  \label{Zcan}
\end{eqnarray}
where the set $\{y_i\}$ is related to the potential via (\ref{Bfactors}). Moreover, 
it holds that
\begin{equation}\label{tZcan}
Z^{can}_{L,N} = \frac{L}{N} y_{1}^{N} \tilde{Z}^{can}_{N,L-N}\left(\{y\}\right)
\end{equation}
where 
\begin{equation}\label{poin} 
\tilde{Z}^{can}_{N,L-N}\left(\{y\}\right)=
\sum_{\bfn=(n_0, \ldots, n_{N-1})\in\tilde{\Omega}^{can}_{N,L-N}} \prod_{i=0}^{N-1} 
\frac{y_{n_{i}+1}}{y_{1}},\,  N \in \N_{1, L},
\end{equation}
 is the canonical partition function of 
zero range process with $L-N$ particles on torus $\T_N$ as defined by (\ref{ZRPcanZ}) with $\bfn=(n_0, \ldots, n_{N-1})$ being the configurations in the canonical state space $\tilde{\Omega}^{can}_{N,L-N}$ for particles without exclusions defined in Section~\ref{SpaceZRP}) and with
the sequence $\{g_i\}$ defining the zero range process
given by
\begin{equation}\label{g-to-y}
g_0=0,\,\, g_k=\frac{y_k}{y_{k+1}}
= \rme^{J(k+1)-J(k)}, \, k=1,2,\ldots .
\end{equation}
\end{prop}
%

Before proving the proposition we note that
the main difficulty lies in the derivation of the identity (\ref{Zcan}). Observe that
it expresses $Z^{can}_{L,N}$ using headways of configurations from $\Omega^{can}_{L,N}$. This  
is clear because of the presence of the constraint $\delta_{\sum_{i=1}^N n_{i},L-N}$ that restricts the summation
to those sequences $(n_1, \ldots, n_{N})$ 
in which the sum is equal to the lattice size $L$ diminished 
by the number of particles $N$ on the lattice, see the relation (\ref{distsum}) that says that this property 
is characteristic for genuine headway sequences. 
Hence, we may informally say that the proof aims 
at transforming the expression for $Z^{can}_{L,N}$ that has been originally written in terms of occupation 
numbers into an expression that uses headways. In a sense, this transformation is similar to the construction
of a Palm representation of a measure. This similarity indicates that the invariance of measure $\mu^{can}_{L,N}$
w.r.t. translations of the torus will be essential for the proof.

Once (\ref{Zcan}) is established, the passage to (\ref{tZcan}) is almost direct. Actually, it is the 
relation (\ref{tZcan}) that is needed below for the calculation of the current-density relation of the AHEP. 

\medskip \proof Trivially, $Z^{can}_{L,0} =1$. Hence we consider in what follows $N>0$.

As we noticed above, the principal aim is to derive the relation (\ref{tZcan}). In that relation, the summation 
is over the space $\tilde{\Omega}^{can}_{N, L-N}$. This space is not isomorphic to $\Omega^{can}_{L,N}$
that is the space of the ``contributors'' to the value of $Z_{L,N}^{can}$ as determined by the definition 
(\ref{Zcandef}). Thus, the principal effort of the whole proof is to split $\Omega^{can}_{L,N}$ into bunches 
so that the bunch set and the space $\tilde{\Omega}^{can}_{N, L-N}$ become isomorphic. The ultimate reason 
why this separation turns out to be possible is the invariance of $\mu^{can}_{L,N}$  w.r.t. translations of $\T_L$.

Due to translation invariance and fixed particle number in the canonical state space with 
$N$ particles one has
\begin{eqnarray}
N Z^{can}_{L,N} 
& = & N \sum_{\bfeta\in\Omega^{can}_{L,N}} \prod_{x=0}^{L-1} 
\prod_{n=0}^{L-N} y_{n+1}^{h_{x,n}(\bfeta)} \nonumber \\         
& = & \sum_{\bfeta\in\Omega^{can}_{L,N}} \partnum(\bfeta) \prod_{x=0}^{L-1} \prod_{n=0}^{L-N} 
y_{n+1}^{h_{x,n}(\bfeta)}  \nonumber\\                                       
& = & \sum_{\bfeta\in\Omega^{can}_{L,N}} \sum_{y\in\T_L} \eta_y \prod_{x=0}^{L-1} \prod_{n=0}^{L-N} 
y_{n+1}^{h_{x,n}(\bfeta)} \nonumber \\                                        
& = & L \sum_{\bfeta\in\Omega^{can}_{L,N}} \eta_0 \prod_{x=0}^{L-1} 
\prod_{n=0}^{L-N} y_{n+1}^{h_{x,n}(\bfeta)}.  \label{Pm4}
\end{eqnarray}
The first line is just the definition of $Z^{can}_{L,N}$. In the second line, we use that $\partnum(\bfeta) = N$
when $\bfeta\in \Omega^{can}_{L,N}$. In the third line, we use the expression \eref{Ndef} of 
$\partnum(\bfeta)$, and in the last line, we use translation invariance of the measure $\mu^{can}_{L,N}$.

We now notice that if $x$ (any lattice site) and $n$ (the headway of the particle that occupies $x$) 
are in the relation $L-n \leq x \leq L-1$,  then, due to (\ref{head3}), 
$\eta_0  y_{n+1}^{h_{x,n}(\bfeta)} = \eta_0$. This fact and (\ref{Pm4}) allow us to continue as follows:
\begin{eqnarray}
Z^{can}_{L,N} & = & \frac{L}{N} \sum_{\bfeta\in\Omega^{can}_{L,N}} \eta_0 
\prod_{x=0}^{L-1} \prod_{n=0}^{L-N} 
y_{n+1}^{h_{x,n}(\bfeta)} \nonumber \\
& = & \frac{L}{N} \sum_{\bfeta\in\Omega^{can}_{L,N}} \left(\prod_{x=0}^{L-1} 
\eta_0 y_{1}^{h_{x,0}(\bfeta)}\right)
\times \left( \prod_{n=1}^{L-N} \prod_{x=0}^{L-1} \eta_0 
y_{n+1}^{h_{x,n}(\bfeta)} \right) \nonumber \\
& = & \frac{L}{N} \sum_{\bfeta\in\Omega^{can}_{L,N}} \left(\prod_{x=0}^{L-1} 
\eta_0 y_{1}^{h_{x,0}(\bfeta)}\right)
\times \left( \prod_{n=1}^{L-N} \prod_{x=0}^{L-n-1} \eta_0 
y_{n+1}^{h_{x,n}(\bfeta)} \right) \nonumber \\
& = & \frac{L}{N} \sum_{\bfeta\in\Omega^{can}_{L,N}} \eta_0 
\prod_{n=0}^{L-N} \prod_{x=0}^{L-n-1} y_{n+1}^{h_{x,n}(\bfeta)} .
\label{Zcan1a}
\end{eqnarray}

{}From now on, we shall use the coordinate presentation of configurations defined in Section~\ref{Subsub:particles}.
In the coordinate presentation, the expression (\ref{Zcan1a}) reads: 
\[ 
Z^{can}_{L,N} =  \frac{L}{N} \sum_{\bfx\in\Omega^{can}_{L,N}} \delta_{x_1,0} 
\prod_{n=0}^{L-N} 
\prod_{x=0}^{L-n-1} y_{n+1}^{h_{x,n}(\bfeta(\bfx))} .
\]  
Above,  $\delta_{x_1,0}$ plays the role of $\eta_0$ in (\ref{Zcan1a}); it annihilates the contribution 
to the sum of all those configurations that have no particle at site $0$. The requirement $x_1=0$ implies 
also that the summation over $\bfx$ reduces to those of them for which $1 \leq x_2 < x_3 < \dots < x_N \leq L-1$.
Using this implication, we write: 
\begin{equation}\label{Zcan2}
Z^{can}_{L,N} 
 =  \frac{L}{N} \sum_{x_2=1}^{L-N+1} \sum_{x_3=x_2+1}^{L-N+2} 
\dots \sum_{x_N=x_{N-1}+1}^{L-1}       
\prod_{n=0}^{L-N}  \prod_{x=0}^{L-n-1} 
y_{n+1}^{h_{x,n}(\bfeta(0,x_2,\dots,x_N))} .
\end{equation}
Next, we consider the headways $n_i\equiv n_i(\bfeta(\bfx)), i\in \N_{1,N}$ (recall from \eref{hwdef}
that $n_i$ means the number of the empty sites between $i$-th particle and the right neighbor particle 
on its right). They satisfy 
\begin{eqnarray*}
h_{x,n}(\bfeta(0,x_2,\dots,x_N))& = &\delta_{x_{2},n+1}\, = \,\delta_{n_{1},n} \delta_{x_{1},x}\mbox{ for } x=0,\\
h_{x,n}(\bfeta(0,x_2,\dots,x_N))& = &\sum_{i=2}^{N-1} \delta_{n_{i},n}\delta_{x_{i},x},\quad 1 \leq x \leq L-n-2, \\
h_{x,n}(\bfeta(0,x_2,\dots,x_N)) & = &  \delta_{x_{N},L-n-1} \, = \,  \delta_{n_{N},n}
 \delta_{x_{N},x} \mbox{ for } x=L-n-1.   
\end{eqnarray*}
The one in the middle just rewrites $n_i$ in the 
terms of the Kronecker-$\delta$. The first one is based upon the obvious fact that $x_2=n_1+1$ when $x_1=0$, 
and the last one upon the fact that $n_N=L-x_N-1$ when $x_1=0$:
Hence 
\[ 
h_{x,n}(\bfeta(0,x_2,\dots,x_N)) = \sum_{i=2}^{N-1} \delta_{n_{i},n} \delta_{x_{i},x},\quad 0 \leq x \leq L-n-1
\]
so that
\[
\prod_{x=0}^{L-n-1} y_{n+1}^{h_{x,n}(\bfeta(0,x_2,\dots,x_N)}  =  y_{n+1}^{\sum_{i=1}^N \delta_{n_{i},n}}
\]
and therefore
\[
\prod_{n=0}^{L-N} \prod_{x=0}^{L-n-1} y_{n+1}^{h_{x,n}(\bfeta(0,x_2,\dots,x_N))} 
= \prod_{i=1}^N y_{1+n_{i}(\bfeta(0,x_2,\dots,x_N)} 
\]
Plugging this into \eref{Zcan2}, using \eref{distsum} 
with $x_1=0$, and shifting summation indices $x_{i} \to m_{i} + x_{i-1}$ in the summation over the particle 
positions leads to
\begin{eqnarray}
Z^{can}_{L,N} 
& = & \frac{L}{N} \sum_{m_2=1}^{L-N+1} \sum_{m_3=1}^{L-N+2-m_2} 
\dots \sum_{m_N=1}^{L-1-m_{N-1}} y_{1+L- \sum_{n=1}^N m_{n}} \prod_{n=2}^N y_{m_{n}}\nonumber \\
& = & \frac{L}{N} \sum_{m_1=1}^{\infty} \sum_{m_2=1}^{\infty} 
\dots \sum_{m_N=1}^{\infty} \left( \delta_{\sum_{i=1}^N m_{i},L} \prod_{j=1}^N y_{m_{j}} \right)\nonumber \\
& = & \frac{L}{N} \sum_{n_1=0}^{\infty} \sum_{n_2=0}^{\infty} 
\dots \sum_{n_N=0}^{\infty}  \left(\delta_{\sum_{i=1}^N n_{i},L-N} \prod_{j=1}^N y_{n_{j}+1} \right) 
\label{Zcan3a}
\end{eqnarray}
The second equality arises because the Kronecker-$\delta$ restricts the infinite sums to the limits in the 
first equation. The third equality results from a shift $m_{i} = n_{i}+1$ of all summation indices. The 
last line in (\ref{Zcan3a}) is identical to (\ref{Zcan}).

The passage from (\ref{Zcan}) to (\ref{tZcan}) follows by simple transformation. First, we note that
the constraint $\delta_{\sum_{i=1}^N n_{i},L-N}$ can be dropped if we restrict the summation over 
all sequences $(n_1, \ldots, n_N)$ of non-negative integers such that $\sum_{i=1}^N n_i=L-N$. Next, we
interchange summation and product and deduce that (\ref{Zcan}) is equal to 
\[
\frac{L}{N} \sum_{\substack{(n_1, \ldots, n_{N}): n_i\geq 0 \forall i \\ \hbox{and }\sum n_i=L-N}}\,\,\,\,
 \prod_{i=1}^{N} y_{n_i}
\]
At the last step, we write $n_{j-1}$ instead for $n_j$ in the above expression. This is necessary because
(due to our construction) $\tilde{\Omega}^{can}_{N,L-N}$ consists of sequences of $N$ non-negative integers
that sum up to $L-N$ but are indexed by the index that runs from $0$ to $N-1$. All this leads to (\ref{tZcan})
and completes to proof.                            \qed  

\subsection{Stationary current}
\label{Sec:current}
This main result of this section, Theorem~\ref{Theo:current} that
expresses particle current as a function of particle density in large volume limit, is based on two pillars: {\em (1)} \ the expression for the stationary particle current at fixed
volume and fixed particle number that is provided by Proposition~\ref{Prop:jstatLN} and {\em (2)} \ properties
of the zero range process, namely, 
{\em 2(a)} \  expectation of the waiting time for jumps of particles in zero range process, and 
{\em 2(b)} \ equivalence of ensembles for zero range process. A brief introduction to the zero range process and the properties {\em 2(a)} and {\em 2(b)}
is provided in Appendix~A. 
We stress that due to the reflection symmetry pointed out in \ref{Rem:reflection} we can state and prove our result
for the TAHEP. Hence in this section, $\mathcal{L}$ is 
the generator of TAHEP. The generalization of the current-density relation to the
AHEP is then an immediate corollary.

We start defining the entities to which our results relate.  
With the functions $\indic_x(\bfeta):=\eta_x$ defined for each $x\in \T_L$ and each $\bfeta=
(\eta_0, \ldots, \eta_{L-1})\in \Omega^{can}_L$ 
we have
\begin{equation}\label{CircRel}
\left(\mathcal{L}\indic_x\right)(\bfeta) = j^{+}_{x-1}(\bfeta) - j^{+}_{x}(\bfeta),\,\,\, x\in \T_L,\,
\bfeta\in \Omega^{can}_L,
\end{equation}
for some function $j^{+}_x(\bfeta), x\in \T_L,\,\bfeta\in \Omega^{can}_L$ which is called {\it instantaneous current} of the TAHEP. This is true because of translation invariance
of the process on the discrete torus and because of
particle number conservation. 
Notice that as a particular consequence of this
relation, if $x=0$ in (\ref{CircRel}), then $x-1=L-1$ as explained in Remark~\ref{Rem:Conven}).
From the action of the generator \eref{AHEPgen} 
we get that 
\begin{eqnarray}
\left(\mathcal{L}\indic_x\right)(\bfeta) & = & \sum_{x'\in\T_L} c_{x'}(\bfeta) (\eta_{x'+1}-\eta_{x'}) 
(\delta_{x,x'} - \delta_{x,x'+1}) \nonumber \\
& = & c_{x-1}(\bfeta) (\eta_{x-1}-\eta_{x}) - c_{x}(\bfeta) (\eta_{x}-\eta_{x+1}), 
\label{GeneratorCurrent02}
\end{eqnarray}
that now allows us to identify the structure of $j_x(\bfeta)$:
\begin{eqnarray}
j^{+}_x(\bfeta) & = & c_{x}(\bfeta)(\eta_{x}-\eta_{x+1}) \label{CurrentF01} \\
& = & \sum_{n=1}^{L-1} w_{n} h_{x,n}(\bfeta) \label{CurrentF02}
\end{eqnarray}
where in the last passage we used the indicator relations \eref{head4} and \eref{head5} 
of Lemma \ref{Lem:head}.  

We consider the stationary current $j^{+}_{L,N}: =\exval{j_x(\bfeta)}^{can}_{L,N}, \, x\in \T_L$, where 
$\exval{\cdot}^{can}_{L,N}$ means -- here and everywhere below  -- the expectations of w.r.t. invariant
measure
$\mu^{can}_{L,N}$ (\ref{Def:CanHead}).  Since this measure 
is translation invariant, $\exval{j_x(\bfeta)}^{can}_{L,N}$ does not depend on $x$ and using \eqref{CurrentF02} we can write
\begin{equation}\label{StCurrent}
j^{+}_{L,N}  =\sum_{n=1}^{L-1} w_{n} \exval{h_{0,n}}^{can}_{L,N} .
\end{equation}

From the stationary current we construct the thermodynamic
limit 
\begin{equation}\label{JlimitVolume} 
j^{+}(\rho) :=  
\lim_{L\to\infty} j^{+}_{L, \floor{\rho L}}, \, \rho\in (0,1)
\end{equation}
of the stationary current
which we call {\em  current-density relation} in order
to emphasize its dependence on the particle density
$\rho$. Nevertheless,
one should keep in mind that $j^{+}_{L,N}$ 
(and therefore $j^{+}(\rho)$)
depends also on the particle jump rates appearing explicitly in 
\eqref{CurrentF02} and implicitly through the 
potential $\{J_k\}$ defining the canonical invariant
measure via \eqref{rates}. 
We point out that the stationary current is a central ingredient in studying hydrodynamic 
limits  \cite{Spoh91,Kipn99,Frit01} which itself, however, is out of the scope of this work.

%
%


\begin{rem}\label{Rem:OutCons}
Regarding the limitation $L\geq 2$ assumed in Proposition~\ref{Prop:jstatLN}, we note that
when $L=1$ then, trivially, $j^{+}_{1,N}=0$ for $N\in\{0,1\}$. 
Also, $j^{+}_{L,L}=0$, i.e., all lattice sites are occupied. These cases are excluded from 
consideration because our final aim is to take the 
thermodynamic limit $L\to\infty$ of  $j^{+}_{L,N}$ when $N=\floor{\rho} L$ 
for $\rho\in (0,1)$.
\end{rem}

\begin{prop}
\label{Prop:jstatLN} 
%
Fix arbitrarily an interaction potential $\{J_k\}$ satisfying (\ref{Jpotential}) and
let the set $\{y_k\}$ be determined by 
$\{J_k\}$ via (\ref{Bfactors}). 
Next, for arbitrary integer  $L\geq 2$ and  
$N\in \N_{1,L-1}$, consider the TAHEP
with $N$ particles on $\T_L$ with
particle jump rates $\{w_n\}$  related to the potential $\{J_k\}$ via (\ref{rates}) with an 
arbitrary $w>0$. 
Then the stationary current (\ref{StCurrent}) is given by
\begin{equation}\label{jLN} 
j^{+}_{L,N} = w \, \frac{N}{L} \,\frac{\tilde{Z}^{can}_{N,L-N-1}\left(\{y\}\right)}
{\tilde{Z}^{can}_{N,L-N}\left(\{y\}\right)} 
\end{equation}
with the canonical partition functions
\[
\tilde{Z}^{can}_{N,K}\left(\{y\}\right), N\in \N, K\in \N
\]
introduced  in (\ref{poin}) in Proposition~\ref{Prop:Zcan}.
\end{prop}
%

\proof First notice that 
\[ 
\sum_{x=0}^{L-1}h_{x,n}(\bfeta) \pi_{L}(\bfeta) = y_{n+1} \frac{\partial}{\partial y_{n+1}} \pi_{L}(\bfeta)
\]
with the Boltzmann weight \eref{BWdef}. Therefore
\[ 
\sum_{x=0}^{L-1} \exval{h_{x,n}}^{can}_{L,N}  = \frac{1}{Z^{can}_{L,N}}\sum_{\bfeta\in\Omega^{can}_{L,N}} 
y_{n+1} \frac{\partial}{\partial y_{n+1}} \pi_{L}(\bfeta) = y_{n+1} \frac{\partial}
{\partial y_{n+1}} \ln{Z^{can}_{L,N}}
\]
which yields
\[
j^{+}_{L,N}  = \frac{1}{L} \sum_{n=1}^{L-1} w_{n} y_{n+1}\frac{\partial}{\partial y_{n+1}} \ln{Z^{can}_{L,N}} 
= \frac{w}{L} \sum_{n=1}^{L-1} y_{n} \frac{\partial}{\partial y_{n+1}} \ln{Z^{can}_{L,N}} 
\]
Elementary algebra involving reshuffling of summation indices in \eref{Zcan3a} and
using $y_{-k}=0$ for $k\in\N_0$ yields after some computation
\[
\sum_{n=1}^{L-1} y_{n} 
\frac{\partial}{\partial y_{n+1}} Z^{can}_{L,N} = \frac{N}{L-1} Z^{can}_{L-1,N}
\]
Now \eref{jLN} follows from Proposition \ref{Prop:Zcan}. \qed

\begin{cor}
The stationary current $j_{L,N}$ of the AHEP is given by
\begin{equation}
j_{L,N} = (r-\ell) j^{+}_{L,N}.
\end{equation}
For the symmetric HEP where $r=\ell$ the stationary current vanishes.
\end{cor}

This corollary is an immediate consequence of the 
definition of the AHEP as a linear superposition of the
TAHEP and the reverse TAHEP with $r=0$ and $\ell=0$ and of the reflection symmetry which implies that the stationary current 
$j^{-}_{L,N}$
of the
reverse TAHEP is given by
$j^{-}_{L,N} = - j^{+}_{L,N}$.

We use Proposition \ref{Prop:jstatLN} to obtain the current-density relation in the thermodynamic limit $L\to\infty$ for the TAHEP evolving 
in its invariant measure.

\begin{theo}
\label{Theo:current} 
Let the interaction potential $\{J_k\}$ satisfying (\ref{Jpotential}) and the corresponding 
set $\{y_k\}$ 
defined by (\ref{Bfactors}) be such that 
\begin{equation}\label{summability}
\left\{{y_n}/{y_{n+1}}, n=0,1,\ldots\right\} \hbox{ is a limited sequence.}
\end{equation} 
Let then $\{y_k\}$ determine the function $F(\cdot)$ via
\begin{equation} \label{fdef}
F(u) = \sum_{k=0}^\infty u^k y_{k+1}, \hbox{ for } u\in (0,u_c)
\end{equation}
where $u_c$ is such that $(0,u_c)$ is the positive part of the domain of convergence of the 
series that defines $F$. Let then  
\[ 
r(u) := u \ddu \ln{F(u)},\,\, u\in (0,u_c)
\]
and define $\phi(r)$ as the inverse of $r(u)$.
Define
\[
\rho_c:= \frac{1}{1+r(u_c)}
\]
that can be $0$ or some value in $(0,1)$ depending upon the values of $\{y_k\}$.
Then the current-particle relation $j^{+}(\cdot)$, defined by (\ref{JlimitVolume}) for the 
TAHEP with particle jump rates $\{w_n\}$ related to the potential $\{J_k\}$ via 
(\ref{rates}) with an arbitrary $w>0$, is given by the relation
\be 
j^{+}(\rho) = \left\{\ba{ll}
\displaystyle w \rho \, \phi\left(1/\rho-1\right) & \rho\in [\rho_c,1) \\[1mm]
\displaystyle w \rho \, \phi(1/\rho_c-1) & \rho\in [0,\rho_c). \ea\right.
\label{jrho}
\ee
\end{theo}

\proof  We make use of properties of zero range process (ZRP) that are presented in 
Appendix~\ref{App:ZRP}. For readers not familiar with the ZRP we recall that both the dynamics of the ZRP and its invariant measures are constructed 
on the basis of the sequence of departure rates $\{g_k\}_0^\infty$ (\ref{gfunc}). 
To match this construction to the TAHEP  we choose
\begin{equation}\label{ZRPjumps}
g_k:=\frac{y_k}{y_{k+1}}=e^{J(k+1)-J(k)},\quad k\in \N_0
\end{equation}
which is the assumption (\ref{g-to-y}) of Proposition~\ref{Prop:Zcan}. Then
$r(\cdot)$ is equal to the strictly
increasing function $\tilde{\rho}(\cdot)$ \eqref{trhodef}
which ensures that the inverse $\phi(\cdot)$
exists.

We start with two integers $N$ and $L$ that we fix arbitrarily but such that $L>N$. Let also  
\begin{equation}\label{DefK}
K:=L-N \in \N.
\end{equation}
We define $N$ 
auxiliary functions $\gamma_0,\gamma_1, \ldots \gamma_{N-1}$ that map $\tilde \Omega^{gc}_N$ 
introduced in (\ref{TiOmegaGC}) to $\R_0^+$ in such a way that $\gamma_i(\bfn)$ gives the particle 
departure rate from site $i\in\T_N$ in the particle configuration $\bfn$ of the zero range process, namely
\begin{equation}\label{Gafun}
\gamma_i(\bfn):= g_{n_i},\quad \bfn=(n_0, n_1, \ldots, n_{N-1})\in \tilde \Omega^{gc}_N, \, i\in \N_{0, N-1}
\end{equation}
This definition makes it clear that the functions $\gamma_i(\bfn)$ as well as the functions  $g_{n_i}$ are random variables 
and we can consider 
their expectation w.r.t.  $\tilde{\mu}^{can}_{N,K}$, the canonical measure \eqref{canens} of the ZRP with $K$ particles
without exclusion
on a torus with $N$ sites. We denote expectations w.r.t.
this measure by $\exval{\exval {\cdot}}^{can}_{N, K}$.
Due to invariance of the measure w.r.t. translations of $\T_N$, the expectation of $\gamma_i$ does not
depend on $i$. Then, by direct substitution, we get that
\begin{equation}\label{promeg}
\exval{\exval {\gamma_i}}^{can}_{N, K}=\exval{\exval {g_{n_0}}}^{can}_{N,K}=
\frac{\tilde{Z}^{can}_{N,L-N-1}}{\tilde{Z}^{can}_{N,L-N}}, \quad \forall i\in\T_N.
\end{equation}
In the first equality, we substituted the random variable $\gamma_0(\bfn)$ by the random variable 
$g_{n_0}$. 
The second equality is technical but direct: We use (\ref{ZRPcanZ})
with the expression (\ref{piDef}) for $\tilde{\pi}_N(\bfn)$. It follows directly from 
(\ref{piDef}) that $g_{n_0}\tilde{\pi}_N(\bfn)$ changes the product $\prod_{m=1}^{n_0}$ in $\tilde{\pi}_N(\bfn)$
to $\prod_{m=1}^{n_0-1}$. This is why $\tilde{Z}^{can}$ in the numerator acquires the index $N, L-N-1$. 
The denominator normalizes $\tilde{\pi}_N$, thus
transforming it into $\tilde{\mu}^{can}_{N, L-N}$. 
Combining (\ref{promeg}) and Proposition~\ref{Prop:jstatLN} allows us to conclude that
\begin{equation}\label{Metag}
j^{+}_{L,N} \equiv j^{+}_{N+K,N} = w \frac{N}{N+K} \exval{\exval{ g_{n_0}} }^{can}_{N,K}.
\end{equation}

Let now $\rho$ be an arbitrary number from $(0,1)$ that we interpret as the density of TAHEP particles.
Set $N=\floor{\rho L}$ and take the thermodynamic limit
$L\rightarrow \infty$ of (\ref{Metag}). By definition, the left-hand side is the stationary current-density relation
$j^{+}(\rho)$ (\ref{JlimitVolume}). 
The right-hand-side of is, in 
terms of the ZRP, equivalent to the limit $N\rightarrow \infty$  when the density of \zparticles\ 
on $\T_N$ is kept at the value $\tilde{\rho}=1/\rho-1$. Due to the equivalence of ensembles described in
Appendix~\ref{Ap-EquiEnsbl}, the projection at the \zparticle\ configuration at site $0$ 
of the measure $\tilde{\mu}^{can}_{N, \tilde{\rho}N}$ converges, as $N\rightarrow \infty$, to the grand canonical measure
$\tilde{\nu}_{\tilde{\zeta}(\tilde{\rho})}$ defined in
\eqref{DefNu}. Since by (\ref{summability}) and by the
choice (\ref{ZRPjumps}) the departure rate $g_k$  is a limited function on $k\in \N_0$, the expectations also converge, i.e., $\exval{\exval{ g_{n_0}} }^{can}_{N,\tilde{\rho}N}$ converges to the expectation of $g_{n_0}$ w.r.t. the measure 
$\tilde{\nu}_{\tilde{\zeta}(\tilde{\rho})}$. The limit expectation is $\tilde{\zeta}(\tilde{\rho})$, as 
it follows from (\ref{DefNu}) by a simple and direct calculation. 
%
%

The conclusion is that $
j^\ast(\rho)=w\rho \tilde{\zeta}(\tilde{\rho})
$
where $\tilde{\zeta}(\tilde{\rho})$ is constructed in (\ref{phirho}). The theorem 
repeats this construction in its assumptions but formulates its assertion through 
$\rho$ that relates to $\tilde{\rho}$ via $1/\rho=1+\tilde{\rho}$ in the range 
$\rho\in[\rho_c,1)$ for which equivalence of ensembles
is proved in \cite{Kipn99}. It follows that
$\tilde{\zeta}(\tilde{\rho}) = \phi(1/\rho-1)$
for $\rho\in[\rho_c,1)$.
For $\rho\in[0,\rho_c)$ \eref{jrho} holds due to equivalence of ensembles 
in the supercritical regime of the ZRP established in \cite{Gros03} so that
$\tilde{\zeta}(\tilde{\rho}) = \phi(1/\rho_c-1)$
for $\rho\in[0,\rho_c)$.
\qed

%
%

\section{Domain random walk}
\label{Sec:duality}
When $\rho_c \neq 0$ one notices in \eref{jrho} a curious behavior:
The current-density relation $j^\ast(\cdot)$ is linear in $\rho$
in the range $[0,\rho_c]$ as if the particles were non-interacting. However, by analogy to the supercritical regime
in the zero range process we conjecture that this 
behaviour does not arise from some effective absence of interactions at low densities
but from
of a coarsening phenomenon \cite{Godr03,Gros03,Godr17} which in the HEP
is associated with phase separation: 
An initial distribution of particles
that is homogeneous on large scale separates over time into empty domains
and domains with critical density. Eventually, only two domains
survive on a large torus, giving rise to a stable phase-separated state
with one region of density zero and one region with density
$\rho_c$. 

Here we provide support for this conjecture
on a microscopic level for the TAHEP defined on $\Z$
with an arbitrary but finite number $N$ of particles.
We prove a reverse duality \cite{Schu23} that expresses the time evolution
of a certain family of initial measures of the $N$-particle TAHEP defined on $\Z$ (that we
call domain measures) in terms of a convex combination of such measures
with weights given by the transition probabilities of the space-reflected $1$-particle TAHEP 
defined on $\Z$ which is simply a totally asymmetric
random walk. On large scales the domain measures concentrate
on a region of critical density which is consistent with
the conjecture of a stable phase separated state and goes beyond it in so far as it suggests that the phase with
critical density moves with a finite average speed rather
than being fixed in space like in other non-equilibrium
phase separation phenomena \cite{Evan98b,Maha20}.

\subsection{Reverse duality}\label{Subsec:RevDuality}
We recall the notion of reverse duality introduced in \cite{Schu23}. 
For a Markov chain $(\omega(t))_{t\in\R_0}$ with state space $\Omega$
consider the transition function
\begin{equation}
p_t(\omega,\zeta) = \bfP^{\omega}(\omega(t)=\zeta)
\label{transfundef}
\end{equation} 
for all $\omega,\zeta\in\Omega$.
The intensity matrix  $Q$ (i.e., the $Q$-matrix in 
the terminology of \cite{Ligg10}) associated with 
$\omega(t)$ has matrix elements
\begin{equation}
Q_{\omega,\zeta} = \ddt p_t(\omega,\zeta)\vert_{t=0}, \quad 
\omega,\zeta\in\Omega
\label{WQ}
\end{equation}
and gives rise to the Kolmogorov backward equation
\begin{equation}
\ddt p_t(\omega,\zeta)  = 
\sum_{\eta\in\Omega} 
Q_{\omega,\eta}  p_t(\eta,\zeta)  
\label{Kbe}
\end{equation}
for all $t\in\R^+$, see \cite{Ligg10} for details.

\begin{df}
\label{Def:Revduality}
Let $\omega(t)$ and $\tilde{\omega}(t)$ be two Markov chains
with countable state spaces $\Omega$ and $\tilde{\Omega}$ 
respectively and
let $R(\cdot,\cdot)$ be a function on
$\tilde{\Omega}\times\Omega$ that is summable over $\tilde{\omega}$ 
for every $\omega\in\Omega$ and summable over $\omega$ for every $\tilde{\omega}\in\tilde{\Omega}$. The processes 
$\omega(t)$ and 
$\tilde{\omega}(t)$ are said 
to be reversely dual to one another with respect to $R$ if for all $t\geq 0$
the relation
\begin{equation}
\sum_{\omega'\in\Omega} R(\tilde{\omega},\omega') 
\bfP^{\omega'}(\omega(t)=\omega) = \sum_{\tilde{\omega}' \in \tilde{\Omega}}
R(\tilde{\omega}',\omega) 
\bfP^{\tilde{\omega}'}(\tilde{\omega}(t)=\tilde{\omega}) . \quad \label{revdualitydef}
\end{equation}
holds for all $\omega\in\Omega$ and all $\tilde{\omega}\in\tilde{\Omega}$. If 
\eref{revdualitydef} holds
then $R(\cdot,\cdot)$ is called duality function for the two processes 
$\omega(t)$ and $\tilde{\omega}(t)$.
\end{df}

We have not specified the range of the function $R$. From
now on we assume for definiteness that $R$ is real valued.

\begin{prop} 
\label{Prop:revD}
Let $Q$ and $\tilde{Q}$ be the intensity matrices of two Markov processes
$\omega(t)$ and $\tilde{\omega}(t)$
with countable state spaces $\Omega$ and $\tilde{\Omega}$ respectively 
and suppose that the intertwining relation
\begin{equation}
R Q = \tilde{Q}^T R
\label{revduality}
\end{equation}
holds with a matrix $R$ with matrix elements $R_{\tilde{\omega},\omega}$,
$\tilde{\omega}\in\tilde{\Omega}$, $\omega\in\Omega$.
Then $\omega(t)$ and $\tilde{\omega}(t)$ are reversely dual with respect
to the duality function $R(\tilde{\omega},\omega) = R_{\tilde{\omega},\omega}$.
\end{prop}

\proof
We define a family of functions 
$r_\omega: \tilde{\Omega}\times [0,\infty) \to \R$ with index $\omega \in \Omega$ 
by
\begin{equation}
(\tilde{\omega},t) \mapsto
r_\omega(\tilde{\omega},t) := \sum_{\omega'\in\Omega} 
R(\tilde{\omega},\omega') 
\bfP^{\omega'}(\omega(t)=\omega).
\end{equation}
This function is the l.h.s. of \eref{revdualitydef}.
Then by \eref{transfundef} 
\begin{equation}
r_\omega(\tilde{\omega},t) := \sum_{\omega'\in\Omega} 
R(\tilde{\omega},\omega') 
p_t(\omega', \omega) = \sum_{\omega'\in\Omega} 
R_{\tilde{\omega},\omega'}
p_t(\omega', \omega)
\end{equation}
and by \eref{Kbe} and \eref{revduality}
\begin{eqnarray}
\ddt r_\omega(\tilde{\omega},t) 
& = & 
\sum_{\omega''in\Omega} 
\sum_{\omega'\in\Omega}  R_{\tilde{\omega},\omega'}
Q_{\omega',\omega''}  p_t(\omega'',\omega) \nonumber \\
& = & 
\sum_{\tilde{\omega}'\in\tilde{\Omega}} \sum_{\omega''in\Omega}
\tilde{Q}_{\tilde{\omega}',\tilde{\omega}} R_{\tilde{\omega}',\omega''}
 p_t(\omega'',\omega)  \nonumber \\
& = & \sum_{\tilde{\omega}'\in\tilde{\Omega}} 
\tilde{Q}_{\tilde{\omega}',\tilde{\omega}}  r_\omega(\tilde{\omega}',t).
\end{eqnarray}
For $t=0$ one has by definition of the transition function
\begin{equation}
r_\omega(\tilde{\omega},0) = 
R(\tilde{\omega},\omega) = R_{\tilde{\omega},\omega}.
\label{initu}
\end{equation}

This is a system of ordinary differential equations for the functions
$r_\omega(\tilde{\omega},\cdot)$.
As pointed out in
\cite[Exercise 2.38]{Ligg10}
the unique solution of this system of differential equations
with initial condition
\eref{initu} is $r_\omega(\tilde{\omega},t) = \sum_{\tilde{\omega}' \in 
\tilde{\Omega}}
R(\omega,\tilde{\omega}') \tilde{p}_t(\tilde{\omega}',\tilde{\omega})$
for each $\omega\in\Omega$ which is the r.h.s. of \eref{revdualitydef}.
Since by definition $r_\omega(\tilde{\omega},t)$
is the l.h.s. of \eref{revdualitydef} the Proposition is proved. \qed

\begin{rem}
(i) In \cite{Schu23} the intertwining relation \eref{revduality} for the intensity matrices 
$Q$ and $\tilde{Q}$ was used as definition of reverse duality. (ii)
The state spaces $\Omega$ and $\tilde{\Omega}$ do not need to have the 
same cardinality so that, in general, the duality matrix $R$ is not a square 
matrix. 
\end{rem}

%

%
%
By construction, the rows of the duality matrix $R$ are indexed by $\tilde{\omega}$ 
and can therefore be interpreted as vectors with components 
$R_{\tilde{\omega},\omega}$ that can be considered as functions $f^{\tilde{\omega}}: \Omega \to \R$
indexed by $\tilde{\omega}\in\tilde{\Omega}$. An interesting 
application of reverse duality arises if each such function
represents a probability measure 
$\mu^{\tilde{\omega}}(\cdot)$ on $\Omega$, i.e., if the matrix element
$R(\tilde{\omega},\omega)=\mu^{\tilde{\omega}}(\omega)$ 
is the probability of a configuration $\omega\in\Omega$.
Reverse duality then provides information about the time evolution
\be
\mu^{\tilde{\omega}}_t  := \mu^{\tilde{\omega}} S_t
\label{mtdef}
\ee
of this measure under the semigroup $S$ associated with the 
process $\omega_t$ as asserted in the following theorem.

\begin{theo}
\label{Theo:measD}
For countable state spaces $\Omega$ and $\tilde{\Omega}$ let
$\mu^{\tilde{\omega}}$ be a family of probability measures on 
$\Omega$ indexed by $\tilde{\omega}\in\tilde{\Omega}$ and let 
$\omega(t)$ and $\tilde{\omega}(t)$ be Markov chains with transition
functions denoted by $p_t(\cdot,\cdot)$ and $\tilde{p}_t(\cdot,\cdot)$
respectively. The following two assertions are equivalent:
\begin{itemize}
\item[(1)] $\omega(t)$ and 
$\tilde{\omega}(t)$ are reversely dual
w.r.t. the duality function 
$R(\tilde{\omega},\omega)=\mu^{\tilde{\omega}}(\omega)$.
\item[(2)] For an initial measure $\mu_0^{\tilde{\omega}}=\mu^{\tilde{\omega}}$ 
the time evolution of this measure under the semigroup $S_t$ is given by
\begin{equation}
\mu^{\tilde{\omega}}_t(\omega) 
= \sum_{\tilde{\omega}'\in\tilde{\Omega}}  
\tilde{p}_t(\tilde{\omega}',\tilde{\omega}) \mu^{\tilde{\omega}'} (\omega)
\end{equation}
for all $\tilde{\omega}\in\tilde{\Omega}$.
\end{itemize}  
\end{theo}

\proof a) Assertion (1) implies assertion (2): This is an immediate consequence of 
\[ 
\mu^{\tilde{\omega}}_t(\omega) 
=(\mu^{\tilde{\omega}} 
S_t)(\omega) =  \sum_{\omega'\in\Omega} R(\tilde{\omega},\omega') 
\bfP^{\omega'}(\omega(t)=\omega)
= \sum_{\omega'\in\Omega} \mu^{\tilde{\omega}}(\omega')
p_t(\omega',\omega)
\] 
and the definition \ref{Def:Revduality} of reverse duality.\\
b) Assertion (2) implies assertion (1): This follows from taking the time derivative and Proposition \ref{Prop:revD}, using also the uniqueness argument in the proof of Proposition \ref{Prop:revD} for the solutions of the resulting system of ordinary differential equations.
 \qed

This time evolution property means that the probability measure 
$\mu^{\tilde{\omega}}_t(\omega)$ at time $t$ of the
process $\omega(t)$ started with an initial measure 
$\mu^{\tilde{\omega}}_0=\mu^{\tilde{\omega}}$ is given by a convex 
combination of measures from
that family $\mu^{\tilde{\omega}}$, $\tilde{\omega}\in\tilde{\Omega}$
with weights given by the transition functions of the dual process.
These transition functions are the probability measures at time $t$ of the
dual process $\tilde{\omega}(t)$ started with initial configuration
$\tilde{\omega}'$. In other words, for a certain family of initial measures
the measure \eqref{mtdef} of the process at any time $t\geq 0$ is determined 
by the transition function of the dual process.

\subsection{Reverse duality for the TAHEP on $\Z$}\label{SecRevDualTAHEP}

To keep track of constants we consider now the TAHEP
with $r=w$ rather than $r=1$ as done above.
The space-reflected $1$-particle TAHEP on $\Z$ starting from $z^\ast\in\Z$ is then
a totally asymmetric simple random walk $x(t)$ with jumps to the left with rate $w$, 
corresponding to the intensity matrix 
\begin{equation}
\tilde{Q}_{x,y} = w( \delta_{y,x-1} - \delta_{y,x})
\label{tQ1def}
\end{equation}
for $x\leq z^\ast$ and $y\leq z^\ast$. 

With this definition we are in a position to state the main result of this section.

\begin{theo}
\label{Theo:PS} 
Let the partition function
\begin{equation}\label{Zdef}
Z :=  \sum_{n=0}^\infty \prod_{k=1}^{n} g_{k}^{-1}
\end{equation}
with $g_k$ given by the jump rates of the TAHEP as in \eqref{ZRPjumps} be 
finite.
Then the TAHEP $\bfx(t)$ with $N$ particles on $\Z$ and starting with the 
leftmost particle on site $x^\ast\in\Z$ is reversely dual to a single totally asymmetric random walk 
$x(t)$ with jumps to the left and starting at $x^\ast$
with respect to the duality function
\begin{equation}\label{dualitymatrix}
R(x,\bfx) = \frac{\delta_{x_1,x}}{Z^{N-1}} \prod_{i=1}^{N-1}
\prod_{k=1}^{x_{i+1}-x_{i}-1}
g_{k}^{-1}
\end{equation}
for $\bfx\in\Omega^\ast_N$ and integers $x\leq x^\ast$.
\end{theo}

\proof Appealing to Proposition \ref{Prop:revD} we prove \eqref{revduality}
for the intensity matrices $Q$ and $\tilde{Q}$ defined by the matrix elements 
\eqref{tQ1def}, \eqref{QNdef} and the duality matrix $R$ defined by the
matrix elements \eqref{dualitymatrix}. This linear algebra proof is computational
and we describe the main steps.\\
(i) With the shorthand
\begin{equation}\label{Cydef} 
C_{\bfy} := \frac{w}{Z^{N-1}} \prod_{i=1}^{N-1} \prod_{k=1}^{y_{i+1}-y_{i}-1} g_{k}^{-1} 
\end{equation}
one obtains from \eqref{tQ1def} and \eqref{dualitymatrix} the matrix elements
\begin{eqnarray*}
(\tilde{Q}^T R)_{x,\bfy}
& = &  \sum_{y\in\Z} \tilde{Q}^T_{x,y} R_{y,\bfy}  \\
& = &  \sum_{y\in\Z} ( \delta_{x+1,y} - \delta_{y,x}) \delta_{y_1,y}
C_{\bfy}  \\
& = &  ( \delta_{x,y_1-1} - \delta_{x,y_1})
C_{\bfy}
\end{eqnarray*}
of the product $\tilde{Q}^T R$ for $\bfy \in\Omega^\ast_N$ and $x\leq x^\ast$.\\
(ii)
To compute the same matrix elements of the product $RQ$
we first note that according to the definitions \eqref{genZnew} and \eqref{ZRPjumps} the 
intensity matrix of the $N$-particle TAHEP on $\Z$ with the choice $w_\infty=w$ has matrix elements
\begin{eqnarray}
Q_{\bfz,\bfy} & = & \sum_{l=1}^{N} w_{n_l(\bfz)}
(\delta_{\bfy,\bfz^{l}} - \delta_{\bfy,\bfz}) \\
& = &
w \sum_{l=1}^{N} g_{z_{l+1}-z_{l}-1}
(\delta_{y_l,z_l+1} - \delta_{y_l,z_l})
\prod_{\stackrel{j=1}{j\neq l}}^{N}  \delta_{y_j,z_j} 
\label{QNdef}
\end{eqnarray}
for all $\bfz\in\Omega^\ast_N$, $\bfy\in\Omega^\ast_N$, and
with the convention $g_{z_{N+1}-p}:=1$ for all $p\in\Z$.
Therefore
\bea
(RQ)_{x,\bfy}
& = & \frac{w}{Z^{N-1}} \sum_{\bfz\in\Omega^\ast_N} \delta_{z_1,x}
\prod_{i=1}^{N-1} \left(\prod_{k=1}^{z_{i+1}-z_{i}-1}
g_{k}^{-1} \right) \nonumber \\
& & \times \sum_{l=1}^{N} g_{z_{l+1}-z_{l}-1}
(\delta_{y_l,z_l+1} - \delta_{y_l,z_l})
\prod_{\stackrel{j=1}{j\neq l}}^{N}  \delta_{y_j,z_j}
\eea
which we split into terms proportional to
\bea 
A_{x,\bfy} & := &  \sum_{l=1}^{N} \sum_{\bfz\in\Omega^\ast_N}
\delta_{z_1,x} g_{z_{l+1}-z_{l}-1}
\prod_{i=1}^{N-1} \left(\prod_{k=1}^{z_{i+1}-z_{i}-1} g_{k}^{-1} \right)
\delta_{y_l-1,z_l}
\prod_{\stackrel{j=1}{j\neq l}}^{N}  \delta_{y_j,z_j} \nonumber \\
B_{x,\bfy} & := &  \sum_{l=1}^{N} \sum_{\bfz\in\Omega^\ast_N}
\delta_{z_1,x} g_{z_{l+1}-z_{l}-1}
\prod_{i=1}^{N-1} \left( \prod_{k=1}^{z_{i+1}-z_{i}-1} g_{k}^{-1} \right)
\prod_{j=1}^{N}  \delta_{y_j,z_j} 
\eea
so that
\[
(RQ)_{x,\bfy} = \frac{w}{Z^{N-1}} \left(A_{x,\bfy} - B_{x,\bfy}\right).
\]
(iii) For the positive contribution $A_{x,\bfy}$ to the matrix element
$(RQ)_{x,\bfy}$ one finds by splitting the sum over $l$
into the term with $l=1$ and the remaining sum over $l$ from 2 to $N$
\bea
A_{x,\bfy} 
& = &
\sum_{\bfz\in\Omega^\ast_N}
\delta_{z_1,x} g_{z_{2}-z_{1}-1}
\prod_{i=1}^{N-1} \left(\prod_{k=1}^{z_{i+1}-z_{i}-1}
g_{k}^{-1} \right)
\delta_{y_1-1,z_1}
\prod_{j=2}^{N}  \delta_{y_j,z_j} 
\nonumber \\ & &
+ \sum_{l=2}^{N}
\sum_{\bfz\in\Omega^\ast_N}
\delta_{z_1,x} g_{z_{l+1}-z_{l}-1}
\prod_{i=1}^{N-1} \left(\prod_{k=1}^{z_{i+1}-z_{i}-1}
g_{k}^{-1} \right)
\delta_{y_l-1,z_l}
\prod_{\stackrel{j=1}{j\neq l}}^{N}  \delta_{y_j,z_j} .
\eea
The product over the indicators given by the Kronecker-$\delta$'s  $\prod_{\stackrel{j=1}{j\neq l}}^{N} 
 \delta_{y_j,z_j}$ allows for replacing the particle positions $z_j$ by the constants $y_j$, except 
 for $j=l$. For $j=l$
one can replace $z_l$ by $y_l-1$. This yields
\bea
A_{x,\bfy} 
& = &
\delta_{y_1-1,x} g_{y_{2}-y_{1}}
 \left[  \left(\prod_{k=1}^{y_{2}-y_{1}}
g_{k}^{-1} \right) 
\prod_{i=2}^{N-1} \left(\prod_{k=1}^{y_{i+1}-y_{i}-1}
g_{k}^{-1} \right) \right]
\nonumber \\ & & 
\times \sum_{\bfz\in\Omega^\ast_N}
 \delta_{y_1-1,z_1} \prod_{j=2}^{N} \delta_{y_j,z_j}
\nonumber \\ & &
+ \sum_{l=2}^{N}
\delta_{y_1,x} g_{y_{l+1}-y_{l}} 
\nonumber \\ & & \times 
 \left[ \left(\prod_{k=1}^{y_{l}-y_{l-1}-2}
g_{k}^{-1} \right)  \left(\prod_{k=1}^{y_{l+1}-y_{l}}
g_{k}^{-1} \right)
\prod_{\stackrel{i=1}{i\neq l-1,l}}^{N-1} \left(\prod_{k=1}^{y_{i+1}-y_{i}-1}
g_{k}^{-1} \right) \right]
\nonumber \\ & & 
\times \sum_{\bfz\in\Omega^\ast_N}
\delta_{y_l-1,z_l} \prod_{\stackrel{j=1}{j\neq l}}^{N}  \delta_{y_j,z_j}.
\eea
As next step we note that the products in the braces $[\dots]$ are 
proportional to $C_{\bfy}$ defined in \eqref{Cydef} with a  factor involving
$g_{\cdot}$, viz.,
\bea 
\left(\prod_{k=1}^{y_{2}-y_{1}}
g_{k}^{-1} \right) 
\prod_{i=2}^{N-1} \left(\prod_{k=1}^{y_{i+1}-y_{i}-1}
g_{k}^{-1} \right) & = & \frac{1}{g_{y_{2}-y_{1}}}  C_{\bfy} \nonumber \\
\left(\prod_{k=1}^{y_{l}-y_{l-1}-2}
g_{k}^{-1} \right)  \left(\prod_{k=1}^{y_{l+1}-y_{l}}
g_{k}^{-1} \right)
\prod_{\stackrel{i=1}{i\neq l-1,l}}^{N-1} \left(\prod_{k=1}^{y_{i+1}-y_{i}-1}
g_{k}^{-1} \right) & = & \frac{g_{y_{l}-y_{l-1}-1}}{g_{y_{l+1}-y_{l}}} C_{\bfy}. 
\eea
Thus $A_{x,\bfy}$ reduces to
\bea
A_{x,\bfy} 
& = &
\delta_{y_1-1,x} C_{\bfy} \, \sum_{\bfz\in\Omega^\ast_N}
\delta_{y_1-1,z_1} \prod_{j=2}^{N}  \delta_{y_j,z_j}
\nonumber \\ & &
+  \delta_{y_1,x} C_{\bfy} \sum_{l=2}^{N}
g_{y_{l}-y_{l-1}-1}
\, \sum_{\bfz\in\Omega^\ast_N}
\delta_{y_l-1,z_l} \prod_{\stackrel{j=1}{j\neq l}}^{N}  \delta_{y_j,z_j}.
\nonumber 
\eea
The only dependence on the particle positions $z_i$ is in the
indicator functions given by the Kronecker-$\delta$'s so that the
sum over $\bfz\in\Omega^\ast_N$ is equal to 1 for all $\bfy\in\Omega^\ast_N$. Shifting the
summation index in the sum over $l$ we arrive at
\[
A_{x,\bfy} 
 = \left(\delta_{y_1-1,x} + \delta_{y_1,x} \sum_{l=1}^{N-1} g_{y_{l+1}-y_{l}-1}\right)\, C_{\bfy}.
\]
(iii) The summation appearing in the negative contribution $B_{x,\bfy}$
to the matrix element
is simpler since the appearance of the product of indicators 
$\prod_{j=1}^{N}  \delta_{y_j,z_j}$ allows for replacing all particle positions 
$z_j$ by $y_j$. Performing the summation over the indicators then immediately yields
\begin{eqnarray*}
B_{x,\bfy} & = &   \sum_{l=1}^{N} \delta_{y_1,x} g_{y_{l+1}-y_{l}-1}\, C_{\bfy}
 \sum_{\bfz\in\Omega^\ast_N} \prod_{j=1}^{N}  \delta_{y_j,z_j} \\
& = &  \sum_{l=1}^{N} \delta_{y_1,x} g_{y_{l+1}-y_{l}-1}\, C_{\bfy}
\end{eqnarray*}
Splitting the summation over $l$ into a sum from 1 to $N-1$ and the remaining term for $N$ leads to 
\[
B_{x,\bfy}  = \sum_{l=1}^{N-1} \delta_{y_1,x} g_{y_{l+1}-y_{l}-1}\, C_{\bfy}
+ \delta_{y_1,x} g_{y_{N+1}-y_{N}-1}\, C_{\bfy}
\]
for all $\bfy\in\Omega^\ast_N$. Since $g_{y_{N+1}-y_{N}-1} = 1$ we arrive at
\[ 
B_{x,\bfy}  = \delta_{y_1,x} \left(1 + \sum_{l=1}^{N-1} g_{y_{l+1}-y_{l}-1}\right)\, C_{\bfy}.
\]
It follows that
\[ 
(RQ)_{x,\bfy} = \frac{w}{Z^{N-1}} 
\left( A_{x,\bfy}  - B_{x,\bfy} \right) = \left(\delta_{y_1-1,x} 
- \delta_{y_1,x} \right)\, C_{\bfy} = (\tilde{Q}^T R)_{x,\bfy}
\]
which concludes the proof of the Theorem. \qed

\subsection{Random walking domain and condensation}\label{WalkingAndCondensation}

To elucidate the significance of this reverse duality we 
introduce a family of measures 
on $\Omega^\ast_N$ which, roughly speaking, define a configuration of $N$ particles
such that there are no particles to the left of $x^\ast$ while to the right of $x^\ast$ 
one has a domain of mean length $N/\rho_c$ with particle density $\rho_c$.

\begin{df}
For finite partition function $Z$ \eqref{Zdef} the measures 
\[
\mu^{x^\ast}(\bfx) := R(x^\ast,\bfx)
\]
on $\Omega^\ast_N$ indexed by $x^\ast\in\Z$ are called domain measures.
\end{df}

More precisely, the domain measures $\mu^{x^\ast}$ give non-zero probability to configurations
where the leftmost particle is at site $x^\ast$ and the subsequent particles are distributed 
at independent distances in the manner such that headway of $i$-th particle attains the value $n_i$ 
with the probability  
\[
p_i = \frac{1}{Z}\prod_{k=1}^{n_i} g_{k}^{-1}.
\]
The partition function $Z$ is convergent for rates of the form
\begin{equation}\label{Zconvergent}
g_k = 1 + \frac{b}{k} + o(1/k)
\end{equation}
with $b>b_c=2$, corresponding to $\rho_c>0$, see \cite{Gros03,Evan05} for a detailed
discussion in the framework of the ZRP. This choice of rates corresponds to
an  interaction potential that approaches a logarithm for large distances and coincides 
with the onset of condensation in the totally asymmetric ZRP for open boundaries \cite{Levi05,Harr06}.

As a consequence of the general time evolution property for reverse duality established 
in Theorem \ref{Theo:measD} and the reverse duality for the TAHEP 
established in Theorem \ref{Theo:PS} the domain measures evolve in time in a particular
simple way as asserted in the following theorem.

\begin{theo}\label{Theo:mixture}
Start the $N$-particle TAHEP 
on $\Z$ with the leftmost particle at $x^\ast$ and particles distributed with initial domain measure
\[
\mu^{x^\ast}_0(\bfx) = \mu^{x^\ast}(\bfx).
\]
Then at all later times $t>0$ the measure of the TAHEP with $N$
particles is the convex combination 
\[  
\mu^{x^\ast}_t(\bfx) = \sum_{y=x^\ast}^\infty \tilde{q}_t(y,x^\ast) \mu^{y} (\bfx)
\]  
of domain measures $\mu^{y}$ where the weight 
\be 
\tilde{q}_t(y,x^\ast) = \rme^{-wt} \frac{(wt)^{y-x^\ast}}{(y-x^\ast)!}.
\label{qxx0}
\ee
is the transition probability of a single totally asymmetric random walk with jumps to the left
with rate $w$ and starting at $y$, or, equivalently, the transition probability of a single totally 
asymmetric random walk with jumps to the right with rate $w$ and starting at $x^\ast$.
\end{theo}

\begin{figure}
\centering
\includegraphics[width=.9\textwidth]{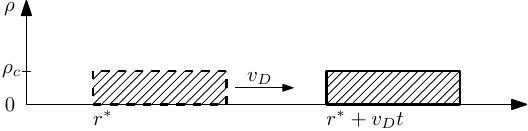}
	\caption{Motion of the domain with density $\rho_c$ on hydrodynamic scale.
The left box marked by broken lines shows the domain at time $t=0$ with the left domain
wall (the shock) at a position $r^\ast$. The right box shows the domain at time $t>0$.
It moves with deterministic velocity $v_D$ without changing its shape. 
} 
	\label{Fig:DRW}
\end{figure}

On large scales, this time evolution property of initial measures from the family $\mu^{y}$ of domain measures 
describes a totally asymmetric random walk of the domain with mean velocity $v_D = w$.
Hence in a single realization of the $N$-particle TAHEP 
started with an initial configuration drawn from the initial domain measure 
$\mu^{\ast}$ one will typically see at a later time a configuration of particles
forming a domain with leftmost particle at $y$, where the random
position $y$ is distributed according to the law of the totally asymmetric
random walk. Thus an initial domain of density $\rho_c$ remains a domain of 
density $\rho_c$, but its position performs a random walk. This phenomenon
occurs for transition rates that lead to condensation in the invariant
measure of the ZRP. On Eulerian hydrodynamic scale
\cite{Reza91,Ferr91} the random walk of the domain
appears as ballistic motion with deterministic velocity $v_D$ as
discussed in the introduction and illustrated in 
Fig. \ref{Fig:DRW}.

\section{Related problems and open questions}
\label{Sec:conc}

Due to its generality the HEP is linked to a very large number of related problems. 
We conclude this work by addressing a few of them and pointing out 
questions that are left open. 

\paragraph{Jump rates:}  If one allows the rates of the AHEP to depend on the length 
of the torus and/or on the total number of particles then the proof of 
Theorem~\ref{Theo:invariance} entails that the uniqueness of the rates 
is lost and a much larger family of rates leaves the Ising 
distance measure with arbitrary interaction potential invariant. 
In fact, the equations that the rates have to satisfy for fixed length $L$
and fixed particle number $N$ are obtained in explicit 
form and may be regarded as interesting in their own right. 
We refrain from a discussion as in the present context we are not interested
in such processes. However, it seems worthwhile to study more general
processes for which the headway measures are invariant for a given
system size and number of particles.

\paragraph{Phase transitions:}
From a small number of lattice gas models that are amenable to mathematically rigorous
or at least well-founded heuristic analysis it is known that long range interactions
in driven systems with local jumps 
may be accompanied by stationary density correlations that decay algebraically 
with distance \cite{Spoh99,Jack15a,Kare17}, in contrast to short range 
interactions that in one-dimensional systems with one conserved density typically lead to short 
range correlations \cite{Garr90}. Long-range correlations 
in driven systems indicate the presence of non-equilibrium phase transitions
which are poorly understood in the sense that there is no clear picture under
which generic conditions such phase transitions may be expected
to occur in one space dimension, except
in the case of nonconservative particle systems with absorbing states
\cite{Henk09} where the voter model \cite{Ligg99} and the equivalent process of 
diffusion-limited annihilation \cite{Lush87,Gryn95,Schu95} 
is a paradigmatic example, or variants of directed percolation \cite{Jans81,Gras82}
of which the contact process is a paradigmatic example \cite{Ligg99}.

The HEP ``imports'' the well-understood condensation transition 
\cite{Bial97,Jeon00,Evan00,Gros03,Evan05,Arme09,Arme13} of the 
zero-range process (which is not accompanied by long-range spatial correlations)
into the spatial structure of the invariant measure and thus implies a phase transition 
leading to a state with long-range correlations. The result on the 
random walking density domain suggests that this happens in the guise of phase 
separation, so far only known to occur in particle systems with more than one conservation law 
\cite{Evan98b,Lahi00,Clin03,Kafr03,Chak16,Maha20}. The exact properties
of the stationary density correlations in the AHEP are not studied here, but the
process
appears to represent a paradigm for nonequilibrium phase transitions
in one-dimensional particle systems with one conservation law.

\paragraph{Hydrodynamic limit:}
The original motivation behind introducing stochastic
lattice gas models with many
identical particles was the derivation of deterministic large scale hydrodynamic behaviour
from the microscopic laws of interaction in which, as mentioned above, the
stationary current-density relation plays a pivotal role.
On general grounds discussed in depth in \cite{Spoh91} and derived rigorously 
for a large 
class of exclusion models in \cite{Kipn99,Frit01} 
one expects in the hydrodynamic limit under Eulerian scaling
the large scale dynamics of the coarse-grained particle density $\rho(x,t)$
to be governed by the non-linear pde
\[                                           
\partial_t \rho(x,t) + 
\partial_x j(\rho(x,t)) = 0
\]                             
where $j^\ast(\cdot)$ is the stationary
current-density relation.

As seen above, in the AHEP there may be a 
non-zero critical density below which the current becomes linear in
particle density. This would yield the linear evolution
equation $\partial_t \rho(x,t) + 
w \partial_x \rho(x,t) = 0$ for {\it arbitrary}
initial density profiles $\rho(x,0)$ which is in conflict
with the coarsening expected from the related
condensation phenomenon of the zero-range 
process (ZRP) where it implies a breakdown of hydrodynamics above some critical density.
In the ZRP the conditions on the rates for condensation to occur are known very precisely. Hence 
we can say equally precisely under which conditions on the interaction parameters $J(k)$ 
condensation (and therefore breakdown of hydrodynamics) below a critical density is expected to occur. 
For the ZRP the precise nature of breakdown is an open problem for supercritical densities \cite{Schu07,Stam15} 
which translates into an similar open problem for the AHEP for subcritical densities. 

\paragraph{Shock discontinuities:} For a convex current-density relation
one expects from the Rankine-Hugoniot condition solutions
with shock discontinuities \cite{Lax73}. In particular,
for the Riemann problem with initial data 
\[           
	\rho_0(x) = \left\{\ba{ll}
	\rho_- & x < x^\ast \\
	\rho_+ & x > x^\ast
	\ea\right.
\]          
one expects for $\rho_+>\rho_-$ the shifted initial profile
\[                    
	\rho(x,t) = \left\{\ba{ll}
	\rho_- & x < x^\ast + v_s t \\
	\rho_+ & x > x^\ast +v_s t
	\ea\right. = \rho_0(x+v_s t)
\]                        
with the shock velocity
\[     
	v_s = \frac{j(\rho_+)-j(\rho_-)}{\rho_+-\rho_-}
\]   
for any range of densities in which $j(\rho)$ is convex.

Indeed, this behavior has been proved rigorously \cite{Kipn99}
for many IPS, the most prominent example being the ASEP 
\cite{Reza91,Ferr91}. The discontinuity at the position $x(t)=x^\ast+v_st$ 
is a stable
shock. On the other hand, for $\rho_+<\rho_-$ (with $j(\rho)$ convex)
such a shock is not stable and the density profile evolves 
into a rarefaction wave. On finer scales, where fluctuations
play a role, the shock position
performs a Brownian motion around its mean position $x(t)$.

Even though the current-density relation of the AHEP is also convex the reverse duality proved in 
Sec.~4                                                    
for a wide class of jump rates suggests that on hydrodynamic scale the Riemann initial condition 
\[                   
	\rho_0(x) = \left\{\ba{ll}
	0 & x < x_- \\
	\rho_c & x_- < x < x_+ \\
	0 & x > x_+
	\ea\right.
\]                   
with the critical density $\rho_c$ allows for the simultaneous existence of
both a stable shock and a stable antishock: Initial data with this density profile
are expected to evolve on hydrodynamic scale into the shifted initial profile
\[                  
	\rho(x,t) = \left\{\ba{ll}
	0 & x < x_- + v_s t\\
	\rho_c & x_- + v_s t< x < x_+  + v_s t\\
	0 & x > x_+ + v_s t
	\ea\right. = \rho_0(x+v_s t)
\]                  
with stable shock and a stable antishock and shock velocity $w$.
This is tantamount to a phase separation
into a domain of density $\rho_c$ and a complementary domain of
density 0, a phenomenon that does not occur in the usual hydrodynamic
limits of particle systems. This result, however, is in agreement
with non-rigorous arguments for the hydrodynamic limit of the supercritical 
ZRP discussed in \cite{Schu07} but for which a rigorous probabilistic derivation
is still an open problem \cite{Stam15}.

\paragraph{Shock and domain random walk:} 
On microscopic scale it was further shown that for specific initial 
densities the shock position performs a biased random walk and
multiple consecutive shocks evolve into a fluctuating
bound state of shocks \cite{Beli02}, which on macroscopic scale
is reflected in coalescence of shocks \cite{Ferr00}. 
As shown in \cite{Schu23} this random walk property
is a consequence of reverse duality. Remarkably,
the phenomenon of random walking shock was
discovered also in other particle systems, see e.g. \cite{Kreb03,Bala10,Bala19}  and has an analogue
also in quantum spin chain where the shock discontinuity
appears in the guise of spin helix states \cite{Popk13,Popk17,Popk21} 
It would be interesting to explore whether also these are examples of reverse or the
closely related intertwining 
duality \cite{Schu23} and whether also domain
random walks can be proved for some SIPS.

\paragraph{Coarsening and metastability:}
The AHEP provides the insight that generally nearest-neighbor interactions 
with logarithmic long-range potential induce phase separation below a 
non-universal critical value. Appealing to the ZRP this can be understood
by noting that when condensation takes place in the ZRP on a finite torus
at density $\rho_c$ then in a measure with a particle density $\rho$ above 
$\rho_c$ all sites but one will be occupied with density $\rho_c$
and one site will carry the excess mass. Correspondingly, in the AHEP
a measure with a density $\rho$ below $\rho_c$ a typical configuration
will have a macroscopic segment on the torus of length 
$L(1-\rho/\rho_c)$ which is completely empty. In complementary segment
of length $L \rho/\rho_c$ particles are distributed
according to the measure at density $\rho_c$.
This insight opens up a series of questions, in
particular how the condensation dynamics of the zero-range process that involves 
metastability \cite{Belt12,Land14,Belt17} and coarsening processes 
\cite{Gros03,Godr03,Godr17,Arme23} can be understood and in a rigorous
fashion for the AHEP. 

\paragraph{Density fluctuations:} For particle systems
with one conservation law it is generally expected that
stationary density fluctuations belong either the to
Edwards-Wilkinson universality class \cite{Krug97}
or to the Kardar-Parisi-Zhang universality class 
\cite{Corw12,Spoh17,Mate21}, the two being the
lowest members of a intriguing discrete family of
dynamical universality classes governed by dynamical
exponents that are the ratios of consecutive Fibonacci
numbers \cite{Popk15b}. It is an open question, however,
what fluctuation pattern exist at the onset of
condensation in the zero range process.  Numerical simulation of an exclusion process with long-range interactions \cite{Priy16} suggest the existence of
a continuously
varying dynamical exponent under certain conditions
on the model parameters.
The AHEP
opens up a new approach to study this problem 
rigorously or numerically which
might be more accessible than the ZRP
due to the bounded local state space.

\paragraph{Biological systems:}
We mention two applications of exclusion processes for the study of processes ongoing in biological cells, viz.,
transport by molecular motors and denaturation of DNA.
For a detailed review and more applications see \cite{Chow24}.\\
(i) Originally, two years before the work by Spitzer, the TASEP was 
introduced as a model for the kinetics of protein synthesis by ribosomes (the 
exclusion particles) moving along an RNA template (the lattice) 
\cite{MacD68} and used to explain an experimentally observed slowing down of 
motion of ribosomes due to a molecular ``traffic jam'', see \cite{Schu97a} for 
a detailed discussion. Over the years the notion of molecular traffic jams, 
analogous to traffic jams in vehicular highway traffic, has gained considerable
traction and led to a vast body of work on the modelling of biological molecular
motors by exclusion processes, see e.g. \cite{Chow24} for a recent overview.
Recently, two of the authors of this work studied an exclusion process 
with next-nearest-neighbour interaction and an internal degree of freedom
which allowed for understanding qualitatively a phenomenon counter acting the 
slowing down due to traffic jams, viz. collective pushing, that leads to an 
enhancement of the average velocity of molecular motors \cite{Beli19} 
that was observed experimentally for RNA polymerase during transcription 
elongation \cite{Epsh03a}. The AHEP studied in the present work may serve as a 
basis for a more quantitative approach to collective pushing as it allows
for adjusting the parameters of the interaction potential to fit experimental data.
In particular, the size of the molecular motors in terms of the size of the
units of the underlying templates \cite{MacD68,Laka03,Shaw03,Scho04} 
can be taken into account.

(ii) Long-range interaction potentials play a role in biological cells also in 
the denaturation process of DNA molecules \cite{Bar09} where a melting transition
has been predicted to occur \cite{Carl02} whose entropic origin has been traced
to an effective logarithmic interaction potential. The relative motion of two 
particle
in such a potential has been shown by scaling analysis of the underlying
Fokker-Planck equation to exhibit a rich scaling behaviour that may depend on 
the initial distribution of the particle distance \cite{Kess10,Hirs11}. 
Hence such anomalous behaviour is expected for two particles in AHEP
with an appropriately chosen logarithmic potential which has recently been 
investigated in terms of a selfsimilar stochastic differential equation in 
\cite{Elia21}. An open problem for future study that may be addressed for the 
AHEP is the description of the large scale properties of a such system of more than
two interacting particles.

\section*{Acknowledgements}
This work is financially supported by CAPES, Finance Code 001, by the FAPESP grant 2023/13453-5,  
by CNPq,  grant number 140797/2018-1, and by FCT/Portugal
through project UIDB/04459/2020 (DOI 10-54499/UIDP/04459/2020, 
and through the grants 2020.03953.CEECIND (DOI 10.54499/2020.03953.CEECIND/CP1587/CT0013, and 2022.09232.PTDC (DOI
10.54499/2022.09232.PTDC). 
N. Ngoc gratefully acknowledges the financial support of CAPES and CNPq during his 
studies at the Doctorate Program in Statistics at the University of S\~ao Paulo.

\appendix

\section{Zero range process (ZRP)}
\label{App:ZRP}

We recall well-known facts about the zero range process 
on the discrete torus that are used in proofs
of the present manuscript.

The zero range process (ZRP) \cite{Spit70,Ligg85,Evan05,Kipn99} is an interacting 
particle system in which identical particles move on lattice sites without the exclusion  
interaction. For the ZRP we shall use the term \zparticles\ in order to 
distinguish them from the AHEP particles. 
For $\T_N$ with $N$ sites, the state space of all \zparticle\ configurations
is the grand canonical state space
$\tilde{\Omega}^{gc}_N$ that has been defined in Section~\ref{SpaceZRP}.
We recall that an element from this space is denoted by 
the $N$-tuple $\bfn =(n_0,\dots,n_{N-1})$ where $n_i$ is interpreted as the 
number of particles that are at the site  $i$ ($i\in \{0,1,\ldots, N-1\}$) of the lattice $\T_N$. 
The dynamics is as follows: To each site $i$, $i=0,1, \ldots, N-1$, we attribute  
a clock whose alarm time is ruled by the exponential distribution with a parameter $g_{n_i}$ that depends 
upon $n_i$, the number of particles at site $i$, but is independent of everything else. 
The numeric sequence 
\begin{equation}\label{gfunc}
g_0=0, \, g_1>0,\, g_2>0, \ldots 
\end{equation}
is considered as a parameter of the process and its members are called {\em departure rates}.  
When the first of $N$ alarms goes off, one of the particles that are at the corresponding site 
jumps instantaneously to another site drawn from some
distribution that we do not need to specify since the invariant measure of the process is independent of that distribution. Immediately after the jump, the rates of the 
exponential distributions are adjusted to the new particle configuration and the whole procedure is repeated. 
%

It is known that the process is well defined provided 
\begin{equation}\label{gassum} 
\sup_{k\in\N_0} |g_{k+1} - g_{k}|<\infty\hbox{ and }g_k > g_0 = 0 \quad \forall \, k > 0. 
\end{equation} 
Hence, a wide range of $g$ is acceptable, however we make use exclusively of the following 
specific choice of $\{g_k\}$: $g_k={y_k}/{y_{k+1}}, \, k\in \N_0$, where $\{y_k\}$ is the sequence 
defined (\ref{Bfactors}), that has been used to construct measure that are invariant for AHEP. The reason
for this choice is explained in the text that precedes  (\ref{ZRPjumps}).
%

\subsection{Invariant measures of the ZRP}\label{Ap-InvarMeas}
The entities related to the ZRP will be marked by $\tilde{}$ to distinguish
between the ``ZRP world'' and the ``AHEP world''. Expectation w.r.t. a measure of the ZRP world 
is denoted $\exval{\exval{\cdot}}$ while if a measure belongs to the AHEP world, the notation
is $\exval{\cdot}$.

To define the invariant measures of the ZRP, we define the numeric sequence
\[
\tilde{\mnogitel}_0:=1,\, \tilde{\mnogitel}_k := \prod_{m=1}^k \frac{1}{g_m}
,\, k=1,2, \ldots,
\]
and the unnormalized product measure
\begin{align}
\tilde{\pi}_N(\bfn) := \prod_{i\in \T_N} \tilde{\mnogitel}_{n_i} \equiv & 
\prod_{i=0}^{N-1} \tilde{\mnogitel}_{n_i}= \prod_{i=0}^{N-1}\left( \prod_{m=1}^{n_i} \frac{1}{g_m}\right),
\label{piDef} \\ 
& N\in \N,\, \bfn=(n_0, \ldots, n_{N-1}) \in \tilde{\Omega}^{gc}_N,\nonumber
\end{align}
with the
convention that products with upper limit $0$ are understood as being equal to 1.
%
For each $N$ and $K$, we normalize $\tilde{\pi}_N(\cdot)$ by 
\begin{equation}\label{ZRPcanZ}
\tilde{Z}^{can}_{N,K} 
:= \sum_{\bfn\in\tilde{\Omega}^{can}_{N,K}} \tilde{\pi}_N(\bfn).
\end{equation}
and define normalized measure 
\begin{equation}\label{canens}
\tilde{\mu}^{can}_{N,K}(\bfn): = \frac{1}{\tilde{Z}^{can}_{N,K}}
\tilde{\pi}_N(\bfn),  
\quad \bfn=(n_0, \ldots, n_{N-1}) \in \tilde{\Omega}^{can}_{N,K}.
\end{equation}
The set $\tilde{\Omega}^{can}_{N,K}$ together with the $\sigma$-algebra and the measure $\tilde{\mu}^{can}_{N,K}(\cdot)$ 
is the {\em canonical ensemble} of the ZRP with $K$ \zparticles\ on $\T_N$. 
The measure (\ref{canens}) is known to be the unique invariant measure for this ZRP but this fact 
is not used directly in our present study. 

The {\em grand canonical ensemble with fugacity $\tilde{z}$} of ZRP on $\T_N$ 
is the probability space
$(\tilde{\Omega}^{gc}_{N}, \tilde{\mathcal{F}}^{gc}_{N}, \tilde{\mu}^{gc}_{N,\tilde{z}} )$
where $\tilde{\mathcal{F}}^{gc}_{N}$ means the $\sigma$-algebra of all subsets 
of ${\tilde{\Omega}^{gc}_{N}}$ and $\tilde{\mu}^{gc}_{N,\tilde{z}}(\cdot)$ is the measure on it 
defined as 
\begin{equation}
\tilde{\mu}^{gc}_{N,\tilde{z}}(\bfn) := \frac{\tilde{z}^{\left|\bfn\right|}  \tilde{\pi}_{N}(\bfn) }
{\tilde{Z}^{gc}_{N,\tilde{z}}} 
\equiv \frac{ \tilde{z}^{(n_0+n_1+\cdots+n_{N-1})}  }{\tilde{Z}^{gc}_{N,\tilde{z}}} 
  \tilde{\pi}_{N}(\bfn), 
\quad \forall\,\bfn\in \tilde{\Omega}^{gc}_{N}, \label{grandcanensdef}
\end{equation}
for all positive $\tilde{z}$ inside the domain of convergence of the {\it grand canonical partition function}
\begin{equation}\label{ZzNormal}
\tilde{Z}^{gc}_{N,\tilde{z}}:= \sum_{\bfn\in \tilde{\Omega}^{gc}_N}\tilde{z}^{\left|\bfn\right|}  
\tilde{\pi}_{N}(\bfn) 
\end{equation}
of the ZRP.
 
Since $\tilde{z}^{\left|\bfn\right|}$ factorizes and since $\tilde{\pi}_N(\cdot)$ factorizes as well, also
 $\tilde{\mu}^{gc}_{N,\tilde{z}}(\cdot)$ is a product
measure. With the normalizing constant  
\begin{equation}
\tilde{\Psi}_{\tilde{z}}:=\sum_{m=0}^\infty \tilde{z}^m \tilde{\mnogitel}_m 
\label{DefPsi}
\end{equation}
and the marginal
\begin{equation}
\tilde{\nu}_{\tilde{z}}(k):= \frac{\tilde{z}^{k} \tilde{\mnogitel}_{k}}{\tilde{\Psi}_{\tilde{z}}}
\label{DefNu}
\end{equation}
the grand canonical measure \eqref{grandcanensdef}
is given by
\begin{equation}\label{marginals}
\tilde{\mu}^{gc}_{N,\tilde{z}}(\bfn) =
\prod_{i=0}^{N-1} \tilde{\nu}_{\tilde{z}}(n_i) ,\quad
\forall\, \bfn=(n_0, \ldots, n_{N-1})\in \tilde{\Omega}^{gc}_{N}
\end{equation}
 inside the domain of convergence of 
$\tilde{\Psi}_{\tilde{z}}$. The factorization implies that
\begin{equation}\label{GCexponent}
\tilde{Z}^{gc}_{N,\tilde{z}} = \left(\tilde{\Psi}_{\tilde{z}}\right)^N.
\end{equation}
A further consequence of the factorization is that the range of $\tilde{z}$ is the non-negative 
part of the domain of convergence of the series that define $\tilde{\Psi}_{\tilde{z}}$. The other 
important consequence is that $\tilde{Z}^{gc}_{N,\tilde{z}}$ and $\tilde{\Psi}_{\tilde{z}}$ are defined 
on the same range of values of $\tilde{z}$.

The relation (\ref{marginals}) suggests that each $\tilde{\mu}^{gc}_{N, \tilde{z}}$ may be viewed as 
a projection on $\tilde{\Omega}^{gc}_N$ of the ``infinite product measure''  
$\tilde{\sigma}_{\tilde{z}}$ that is defined on $\left(\N_0\right)^\infty$ as a product measure 
which marginals are distributed as $\tilde{\nu}_{\tilde{z}}$. Specifically,  
\begin{equation}\label{ProdM}
\tilde{\sigma}_{\tilde{z}}\left(\mathcal{A}\right)=\prod_{j=0}^m\tilde{\nu}_{\tilde{z}}(k_j),
\end{equation}
for any $m$, any $k_0, \ldots, k_m$, and event $\mathcal{A}:=\{\bfn=(n_0, n_1, \ldots)\in 
\left(\N_0\right)^\infty\, \hbox{ such that }n_0=k_0, \ldots, n_m=k_m\}$. Notice that also $\tilde{\sigma}_{\tilde{z}}$ is well-defined for  $\tilde{z}$ from the positive 
part of the domain of convergence of $\tilde{\Psi}_{\tilde{z}}$.

\subsection{Particle density in the grandcanonical ensemble}\label{Ap-Density}

The particle density
of the ZRP in the grand canonical measure,
denoted by $\tilde{\rho}(\tilde{z})$, is the expectation w.r.t. $\tilde{\mu}^{gc}_{N,\tilde{z}}$ of the number of \zparticles\ on $\T_N$ 
divided by the lattice size $N$. Due to the factorization
of the grand canonical measure this quantity does not
depend on $N$ and is given by 
\begin{equation}
\tilde{\rho}(\tilde{z})= \frac{1}{N}\exval{\exval{n_0+\ldots n_{N-1}}}_{ \tilde{\mu}^{gc}_{N,\tilde{z}}}=
\exval{\exval{n_0}}_{\tilde{\nu}_{\tilde{z}}}
= \frac{1}{\tilde{\Psi}_{\tilde{z}}}\sum_{k=0}^\infty k \, \tilde{z}^k\mnogitel_k 
= \tilde{z} \frac{\rmd }{\rmd \tilde{z}} \ln{\tilde{\Psi}_{\tilde{z}}},\label{trhodef}
\end{equation}
for $\tilde{z}\in (0, \tilde{z}_c)$, i.e., in the positive part of the radius of convergence of $\tilde{\Psi}_{\tilde{z}}$.

To discuss the density as a function of the fugacity $\tilde{z}$ and the
role of the radius of convergence we draw on \cite{Kipn99} and
\cite{Gros03}. The range of the density $\tilde{\rho}(\tilde{z})$ as
a function of the fugacity is the interval 
$[0,\tilde{\rho}_c )$, with $\tilde{\rho}(0)=0$ and {\em critical density} 
$\tilde{\rho}_c :=\lim_{\tilde{z}\nearrow\tilde{z}_c} 
\tilde{\rho}(\tilde{z})$. The density $\tilde{\rho}(\tilde{z})$ is strictly increasing 
in its argument $\tilde{z}$ so that the inverse function with domain  $[0,\tilde{\rho}_c )$ which we 
denote by $\tilde{z}(\tilde{\rho})$ is well-defined. Lemma~2.3.3
in \cite{Kipn99} asserts that if $\tilde{z}_c =\infty$, 
then $\tilde{\rho}_c =\infty$, whereas for $\tilde{z}_c <\infty$, both $\tilde{\rho}_c =\infty$ 
and $\tilde{\rho}_c <\infty$ are possible. In the second case, $\tilde{\Psi}_{\tilde{z}_c}<\infty$ 
which implies that the measures $\tilde{\mu}^{gc}_{N, \tilde{z}_c}$, $\tilde{\nu}_{\tilde{z}_c}$ and 
$\tilde{\sigma}_{\tilde{z}_c}$ are well defined, and also that 
 $\tilde{\rho}(\tilde{z}_c) =\tilde{\rho}_c$, i.e., the continuity of $\tilde{\rho}(\cdot)$ at 
 $\tilde{z}_c$.

These facts allow us to the function define $\tilde{\zeta}: [0, \infty)\rightarrow 
[0, \infty)$ as follows
\begin{equation}\label{phirho}
\tilde{\zeta}(\tilde{\rho}):=\left\{\begin{array}{ll}\mbox{inverse of }\tilde{\rho}(\tilde{z})\hbox{ that 
is, } \tilde{z}(\tilde{\rho}) &,
\quad\mbox{for }\tilde{\rho} <\tilde{\rho}_c\\ \qquad \tilde{z}_c &,\quad\mbox{for }
\tilde{\rho}\geq\tilde{\rho}_c \ea\right.\ ,
\end{equation}
In this definition, $\tilde{\rho}_c$ can be either infinite or finite; in the latter case, 
$\tilde{\zeta}(\cdot)$ is continuous at the critical density.  

\subsection{Equivalence of ensembles}\label{Ap-EquiEnsbl}
In the limit $N\to\infty$ with fixed \zparticle\ density, i.e., in the limit 
$N,K\to\infty$, where $\tilde{\rho}=K/N$ is kept fixed, 
the canonical measure (\ref{canens}) is equivalent to the grand canonical measure with
fugacity $\tilde{\zeta}\left(K/N\right)$. 
This equivalence is given through the pointwise limit of the $m$-point marginal defined as
\begin{eqnarray*}
\tilde{\mu}^{can}_{N,K} \left[ \mathbf{x}_m,\mathbf{k}_m\right] &:=&\\
&&\mbox{\kern-8em}\tilde{\mu}^{can}_{N,K}\left(\bfn=(n_0, n_1, \ldots, n_{N-1})\in \tilde{\Omega}^{can}_{N,K}
\, :\, 
n_{x_1}=k_1 ,\ldots ,n_{x_m}=k_m \right). 
\end{eqnarray*}
Equivalence of measures then means that for any $\tilde{\rho}>0$, any $m\in \N$, any vector 
$\mathbf{x}_m$ that fixes marginals at the sites $x_1,x_2, \ldots, x_m$, and any vector $\mathbf{k}_m$
that fixes the number of \zparticles\ at those sites to be, respectively, $k_1, k_2,\ldots, k_m$, it holds that 
\[
\lim\limits_{\substack{{N,K\to\infty}\\ {K/N=\tilde{\rho}}}} 
\tilde{\mu}^{can}_{N,K} \left[\mathbf{x}_m,\mathbf{k}_m\right]=\tilde{\sigma}_{\tilde\zeta(\tilde{\rho})}
\left[\mathbf{x}_m,\mathbf{k}_m\right] \equiv 
\prod_{i=1}^m \tilde{\nu}_{\tilde{\zeta}(\tilde{\rho})} (k_i )
\]
A rigorous result on this equivalence is available in \cite{Kipn99}, Appendix 2, but it does 
not cover the supercritical case when $\tilde{\rho}_c <\infty$ and $\tilde{\rho}\geq\tilde{\rho}_c$. 
The supercritical case is covered in \cite{Gros03}. 

Equivalence of ensembles in the supercritical case may look counter intuitive. 
It can be understood by noting that for a particle density $\tilde{\rho}$ above 
$\tilde{\rho}_c$ the phenomenon of {\it condensation} may occur, see \cite{Evan05} 
for a review from a statistical physics perspective and the recent work
\cite{Arme23} and references therein for a probabilistic discussion. 
Condensation means that for a typical configuration has 
a critical background density $\tilde{\rho}_c$ on all but one site and
the excess particles condense on a single random site \cite{Gros03}
and form the so-called condensate. In an infinite system, however, the 
probability of having this site in any finite subset of sites is zero so that 
the condensate becomes ``invisible'' to expectations of cylinder functions
which only ``feel'' the critical background which is distributed according
the to well-defined product measure at the critical density.

\end{document}